\documentclass[smallcondensed,envcountsect,envcountsame]{svjour3}
\smartqed
\journalname{Mathematical Programming}
\renewcommand{\makeheadbox}{}

\usepackage{amsmath,amssymb,amsfonts}
\usepackage{graphicx}
\usepackage{booktabs}
\usepackage{tabularx}
\usepackage{xcolor}
\usepackage{algorithm}
\usepackage{algorithmicx}
\usepackage{algpseudocode}
\usepackage{mathtools}
\usepackage{microtype}
\usepackage{etoolbox}
\usepackage{hyperref}
\spnewtheorem{assumption}[theorem]{Assumption}{\bfseries}{\rmfamily}
\AtEndEnvironment{proof}{\qed}

\newcommand{\R}{\mathbb{R}}
\newcommand{\Fc}{\mathcal{F}}
\newcommand{\ip}[2]{\langle #1,#2\rangle}
\newcommand{\norm}[1]{\lVert #1\rVert}
\DeclareMathOperator{\spn}{span}
\DeclareMathOperator*{\argmin}{arg\,min}

\makeatletter\let\@font@warning\@font@info\makeatother
\hypersetup{bookmarksdepth=2,hypertexnames=false,
  pdftitle={Forgetting matters: when subspace optimization accelerates without the strong convexity parameter},
  pdfauthor={Ilya Kuruzov, Ivan Dorofeev, Alexander Gasnikov}}
\begin{document}

\title{Forgetting matters: when subspace optimization accelerates without the strong convexity parameter}
\titlerunning{Forgetting matters}
\author{Ilya~Kuruzov \and Ivan~Dorofeev \and Alexander~Gasnikov}
\authorrunning{I. Kuruzov et al.}
\institute{I. Kuruzov (corresponding author), A. Gasnikov \at
  Moscow Independent Research Institute of Artificial Intelligence, Moscow, Russia \\
  Innopolis University, Innopolis, Russia \\
  \email{kuruzov.ia@phystech.edu}
  \and
  I. Dorofeev \at Independent researcher, Moscow, Russia}
\date{}
\maketitle

\begin{abstract}
Sequential subspace optimization (SESOP) needs no problem constants, is accelerated on smooth convex functions and often accelerates in practice. One of our main contributions is to disprove the conjecture that it is uniformly accelerated on strongly convex functions. We develop a taxonomy of memory mechanisms by window length and forgetting rate. We consider exact subspace methods whose memory is a window of recent gradients and steps plus finitely many fixed cumulative sums. For this class, uniform acceleration requires either a window of length at least of order the square root of the condition number or a sum that forgets at a rate at least of order the inverse square root. For one fixed-factor filter pair with a bounded window this forgetting scale is necessary and sufficient up to constants. SESOP remembers too much and the plane form of conjugate gradients too little: neither is uniformly accelerated, even on infinitely differentiable functions quadratic near the minimizer. A short window converges locally quadratically over blocks when the Hessian is Lipschitz near the minimizer and equals a multiple of the identity plus a low-rank term there. Constructive upper bounds reach this forgetting scale without problem constants. Geometric filter pairs created at fixed times give uniform acceleration with logarithmic memory, also with inexact subproblems. A single pair with a smoothly adapted factor uses two vectors and attains an accelerated linear rate up to a prefactor polynomial in the condition number.
\end{abstract}

\section{Introduction}\label{sec:intro}

We consider $\min_{x\in\R^n}f(x)$ for a $\mu$-strongly convex $f$ with $L$-Lipschitz gradient and condition number $\kappa=L/\mu$. Optimal methods need of order $\sqrt\kappa\log(1/\varepsilon)$ gradient evaluations to reduce the gap by the factor $\varepsilon$, and gradient descent needs $\kappa\log(1/\varepsilon)$ \cite{nemirovsky1983problem,nesterov2018lectures}. For $\kappa=10^6$ this is a factor of a thousand. Nesterov's method \cite{nesterov1983method,nesterov2018lectures} attains the optimal order, but its momentum depends on $\mu$. Unlike $L$, the constant $\mu$ cannot be found by backtracking: finitely many values and gradients never certify a positive lower bound on the curvature. The usual remedy is restarting \cite{odonoghue2015adaptive,roulet2020sharpness,renegar2022simple}.

Conjugate gradients avoid $\mu$ on quadratics. They minimize $f$ exactly over the plane of the gradient and the previous step, choosing the coefficients from the problem data without knowledge of $\mu$ or $L$ \cite{hestenes1952methods}. Exact subspace methods $x_{k+1}=\argmin\{f(x):x\in x_k+S_k\}$ extend this idea. Nemirovski's orth-method \cite{nemirovski1982orth}, known as SESOP \cite{narkiss2005sequential}, uses
\[
S_k=\spn\Bigl\{\nabla f(x_k),\ x_k-x_{k-1},\ x_k-x_0,\ \sum_{i=0}^{k}w_i\nabla f(x_i)\Bigr\}
\]
and attains the rate $O(LR^2/k^2)$ on smooth convex functions without any parameter. It performs well on strongly convex problems \cite{guminov2023accelerated}, in our experiments far better than Nesterov's method with true constants. To the best of our knowledge, it is open whether unrestarted SESOP is accelerated \emph{uniformly} over $\Fc_{\mu,L}$, meaning that
\[
f(x_k)-f^\star\le C\exp\bigl(-ck/\sqrt\kappa\bigr)(f(x_0)-f^\star)
\]
for absolute constants $C,c>0$, every condition number, every function in the class, every starting point and every $k$, where $f^\star=\min_x f(x)$. In this paper we continue the study of subspace methods started in \cite{kuruzov2024some}. Oracle lower bounds do not separate SESOP from accelerated methods, and on quadratics it coincides with conjugate gradients, so a hard instance must carry an infinite trajectory consistent with SESOP's anchor and aggregate. For nonlinear conjugate gradients worst-case limitations are known \cite{nemirovsky1983problem,karimi2012cgso,dasgupta2025nonlinear}.

Our approach starts from the fact that exact minimization over $x_k+S_k$ inherits the rate of an accelerated method as soon as $S_k$ contains its momentum. This momentum is a geometric average of past gradients and steps with a horizon of order $\sqrt\kappa$, so, roughly speaking, acceleration without $\mu$ comes down to memory: which subspaces built from the history contain such an average when $\kappa$ is unknown? SESOP keeps everything, stale directions included, while a memory that forgets too fast cannot look far enough back. For one fixed filter pair with a bounded window, we identify the necessary forgetting scale and construct accelerated methods on that scale. We then reach this scale without knowledge of $\kappa$.

\subsection{Contributions}

\begin{enumerate}
\item \emph{Remembering too much} (Section~\ref{sec:slow}). We measure a memory fixed in advance by its window $m$ and the forgetting index $\delta$ of its cumulative sums. For $m\ge1$ and slowly forgetting sums, whenever $\kappa\ge C_0\max\{8,m+2\}^2$ for a sufficiently large absolute constant $C_0$, we construct functions in $\Fc_{\mu,L}$ that require $\Omega((\kappa/m)\log(1/\varepsilon))$ iterations as $\varepsilon\to0$ for every exact subspace method with this memory (Theorem~\ref{thm:slow}). SESOP and the plane form of conjugate gradients are such methods and need $\Omega(\kappa\log(1/\varepsilon))$ iterations (Corollary~\ref{cor:sesop}).
\item \emph{Forgetting too fast and the right time scale} (Sections~\ref{sec:fastpair} to~\ref{sec:smooth}). A filter pair that forgets much faster than at rate $\kappa^{-1/2}$ leaves a method with a bounded window as slow as a short window (Theorem~\ref{thm:fast}). So uniform acceleration needs a window of length $\Omega(\sqrt\kappa)$ or $\delta=\Omega(\kappa^{-1/2})$, and with a bounded window and one pair exactly $1-\rho=\Theta(\kappa^{-1/2})$ (Corollary~\ref{cor:memA}). SESOP misses this time scale. Like the plane form, it stays slow with one aggregate on $C^\infty$ functions that are quadratic near the minimizer (Propositions~\ref{prop:plane3}, \ref{prop:smooth}).
\item \emph{Windows} (Section~\ref{sec:window}). A window of length $1.39\sqrt\kappa$ is accelerated, so the bound $\Omega(\sqrt\kappa)$ is sharp (Proposition~\ref{prop:window}). A short window converges quadratically over blocks near a minimizer where the Lipschitz Hessian is a multiple of the identity plus a low-rank term (Theorem~\ref{thm:window-ub}). Near such a minimizer this is faster than any linear rate.
\item \emph{Reaching the time scale without $\kappa$} (Sections~\ref{sec:filters} and~\ref{sec:inexact}). Pairs of geometric filters store a horizon of order $\sqrt\kappa$ in two vectors (Lemma~\ref{lem:pair}). Created at times $2^j$, they give a uniformly accelerated method with memory $O(\log k)$ (Theorem~\ref{thm:schedule}). One pair with a smoothly adapted horizon and no resets is accelerated up to a polynomial prefactor, and with the SESOP directions it keeps the convex rate and is conjugate gradients on quadratics (Theorem~\ref{thm:afb}). The schedule also allows inexact subproblems (Proposition~\ref{thm:inexact}).
\end{enumerate}

The global accelerated bounds for windows and filters follow from one potential certificate. Restart schemes are nearly optimal \cite{renegar2022simple,wu2026parameter}, and our rates are no better, but our filter methods need neither $\mu$ nor $L$, are monotone and never reset their memory (Table~\ref{tab:methods}).

The rest of the paper is organized as follows. Sections~\ref{sec:related} and~\ref{sec:prelim} cover related work and the setting, Sections~\ref{sec:main} and~\ref{sec:upper} the lower and upper bounds, and Section~\ref{sec:experiments} experiments.

\begin{table}[t]
\caption{Methods of this paper. Their iteration formulas use neither $\mu$ nor $L$; choosing an accelerated window length requires a bound on $\kappa$, whereas the activation schedule and smooth pair require none. Storage counts vectors besides $x_k$ and $\nabla f(x_k)$, and a filter with factor $\rho$ has horizon $1/(1-\rho)$. ``Up to poly$(\kappa)$'': an accelerated rate with a prefactor polynomial in $\kappa$. The smooth-pair horizon bound is stated for $\kappa\ge4$.}\label{tab:methods}
\centering\footnotesize\setlength{\tabcolsep}{4pt}
\begin{tabularx}{\textwidth}{@{}llll>{\raggedright\arraybackslash}X@{}}
\toprule
method & storage & memory, horizon & weights & accelerated\\
\midrule
SESOP & $3$ & slow, $\infty$ & fixed & no (Cor.~\ref{cor:sesop}), $O(1/k^2)$ if $\mu=0$\\
window $m$ & $2m-1$ & window, $m$ & fixed & iff $m\gtrsim\sqrt\kappa$ (Cor.~\ref{cor:memA}, Prop.~\ref{prop:window})\\
activation schedule & $2\lfloor\log_2k\rfloor$ & fast, $\le k$ & fixed times & uniform (Thm.~\ref{thm:schedule})\\
smooth pair & $2$ & fast, $O(\sqrt\kappa)$ & adaptive & up to poly$(\kappa)$ (Thm.~\ref{thm:afb})\\
smooth pair + SESOP & $5$ & fast and slow & adaptive & same, $O(1/k^2)$ if $\mu=0$ (Thm.~\ref{thm:afb})\\
\bottomrule
\end{tabularx}
\end{table}

\section{Related work}\label{sec:related}

Exact minimization over low-dimensional subspaces is as old as conjugate gradients \cite{hestenes1952methods} and their nonlinear extensions \cite{hager2006survey}. A subspace containing the momentum of an accelerated method \cite{nesterov2018lectures} gives an exact step at least as good as the accelerated one, which underlies the orth-method \cite{nemirovski1982orth}, methods with a small-dimensional relaxation oracle \cite{nesterov2021primal,guminov2019accelerated} and the span search of Drori and Taylor \cite{drori2020efficient}. Narkiss and Zibulevsky \cite{narkiss2005sequential} named the approach SESOP and solved the subproblems by Newton's method, so that for $\varphi(Ax)$ a step costs as many products with $A$ as a gradient step \cite{zibulevsky2010l1}. SESOP has been extended to weakly quasi-convex functions \cite{guminov2023accelerated}, inexact gradients and restarted conjugate gradients under quadratic growth \cite{kuruzov2021sequential,kuruzov2024some}. CGSO \cite{karimi2012cgso} already attains $O(\sqrt\kappa\log(1/\varepsilon))$ without constants of the function with a subspace of dimension $O(\log k)$, using progress tests and anchor and gradient aggregates indexed by dyadic blocks. Its restarts change the origin of the progress tests without changing the iterates. Our contribution combines lower bounds for fixed memories with an activation schedule that uses no progress test and a single adaptive filter pair; the pair has constant storage and an accelerated rate up to a polynomial prefactor. Other methods alternate conjugate gradient and accelerated steps \cite{karimi2024nonlinear} or minimize a model built from past gradients \cite{nesterov2022gradient}. The limited-memory common-directions method \cite{lee2022limited} is accelerated when the subspace contains Nesterov's momentum, which requires $\mu$ and $L$, whereas our filter pairs contain such a momentum without the constants. Tupitsa \cite{tupitsa2026cacute} replaces the exact subproblem of SESOP by a cubic model and attains a local accelerated rate with adaptive restarts.

Restarting a convex accelerated method with period of order $\sqrt\kappa$ gives an accelerated linear rate \cite{nemirovskii1985optimal}. Adaptive schemes decrease an estimate of $\mu$ when progress is insufficient \cite{nesterov2013gradient,lin2015adaptive}, test the progress directly \cite{odonoghue2015adaptive} or search over a grid of constants \cite{roulet2020sharpness,renegar2022simple}. Each restart throws the momentum away. The online estimate of $\mu$ in \cite{malitsky2020adaptive} has no guarantee. The survey \cite{daspremont2021acceleration} listed the estimation of $\mu$ without restarts as an open problem, which NAG-free \cite{cavalcanti2025adaptive} later addressed with a global rate of gradient descent type and acceleration near the minimizer. Our smooth pair estimates no constant: it adapts the forgetting horizon of the memory and is globally accelerated, up to a prefactor polynomial in $\kappa$. Lan, Ouyang and Zhang \cite{lan2026optimal} reach $\norm{\nabla f}\le\varepsilon$ with $O(\sqrt\kappa\,(\log(\mu_0/\mu)+\log(\norm{\nabla f(x_0)}/\varepsilon)))$ gradient evaluations from a guess $\mu_0\ge\mu$. These methods count gradient evaluations, whereas we count subspace steps and bound the memory from below.

Oracle lower bounds \cite{nemirovsky1983problem,nesterov2018lectures} hold for all first-order methods and are attained when $\mu$ is known. Nemirovski \cite{nemirovski1989complexity} extended them to methods that query the oracle on an affine manifold of fixed dimension at every step, covering exact subspace steps. In large dimension they still need of order $\sqrt\kappa\log(1/\varepsilon)$ steps on $\Fc_{\mu,L}$, so acceleration is worst-case optimal in this fixed-dimensional subspace-oracle model. The bounds of Section~\ref{sec:main} hold in a fixed dimension as $\varepsilon\to0$ and decide which memories reach it. Oblivious methods, whose coefficients do not depend on the data, cannot be accelerated without $\mu$ \cite{arjevani2016iteration}. Exact minimization is not oblivious, and our bounds cover it for memories fixed in advance. They rest on the interpolation theorem of Taylor, Hendrickx and Glineur \cite{taylor2017smooth}, which underlies performance estimation \cite{drori2014performance} and the exact rate of GD\nobreakdash-ELS \cite{deklerk2017worst}, applied to infinite trajectories whose residuals keep a uniform margin (Lemma~\ref{lem:bridge}).

\section{Setting and two tools}\label{sec:prelim}

Inner products and norms are Euclidean. The class $\Fc_{\mu,L}(\R^n)$ consists of differentiable $\mu$-strongly convex functions with $L$-Lipschitz gradients, $0<\mu\le L$, and $\kappa=L/\mu$, where the constants are valid bounds rather than the best ones. For a method generating points $x_k$ we write $s_k=x_k-x_{k-1}$ for $k\ge1$ and $s_0=0$, and $f^\star$ and $x^\star$ are the minimum value and the minimizer. For $0<\varepsilon<1$ the complexity $N(\varepsilon)$ is the first $k$ with $f(x_k)-f^\star\le\varepsilon(f(x_0)-f^\star)$. All sets of historical indices are restricted to nonnegative integers.

Unless stated otherwise, the methods perform the exact subspace step
\begin{equation}\label{eq:step}
x_{k+1}=\argmin_{x\in x_k+S_k}f(x),\qquad \nabla f(x_k)\in S_k,
\end{equation}
where the linear subspace $S_k$ is built from the history by some rule. One iteration means one exact solution of \eqref{eq:step}, whatever the inner solver. The window method of length $m\ge1$ searches over the last $m$ gradients and steps, $S_k=\spn\{\nabla f(x_i),s_i:k-m<i\le k\}$. GD\nobreakdash-ELS, gradient descent with exact line search, uses $S_k=\spn\{\nabla f(x_k)\}$, and the plane form of conjugate gradients uses $S_k=\spn\{\nabla f(x_k),s_k\}$. SESOP with weights $w_1,\dots,w_d$ uses $S_k=\spn\{\nabla f(x_k),s_k,x_k-x_0,\sigma^1_k,\dots,\sigma^d_k\}$ with the aggregates $\sigma^a_k=\sum_{i=0}^{k}w_a(i)\nabla f(x_i)$. Common choices are Nemirovski's weights, $w(k)=\frac12+\sqrt{\frac14+w(k-1)^2}$ with $w(0)\in\{0,1\}$, and the linear weights $w(i)=i$.

Our upper bounds use only the following properties of \eqref{eq:step}.

\begin{lemma}[exact subspace descent]\label{lem:descent}
Let $f\in\Fc_{\mu,L}(\R^n)$ and let $x_{k+1}$ be given by \eqref{eq:step}. Then $x_{k+1}$ is well defined, $s_{k+1}\in S_k$, $\nabla f(x_{k+1})\perp S_k$, $f(x_{k+1})\le f(x_k)-\norm{\nabla f(x_k)}^2/(2L)$, $f(x_{k+1})-f^\star\le(1-\kappa^{-1})(f(x_k)-f^\star)$ and $2\mu(f(x_k)-f^\star)\le\norm{\nabla f(x_k)}^2\le2L(f(x_k)-f^\star)$.
\end{lemma}
\begin{proof}
Strong convexity gives a unique minimizer $x_{k+1}$, so $s_{k+1}\in S_k$, and optimality on $x_k+S_k$ gives $\nabla f(x_{k+1})\perp S_k$. Since $x_k-\nabla f(x_k)/L\in x_k+S_k$, the descent lemma gives $f(x_{k+1})\le f(x_k-\nabla f(x_k)/L)\le f(x_k)-\norm{\nabla f(x_k)}^2/(2L)$. As $f^\star\le f(x_{k+1})$, this also yields $\norm{\nabla f(x_k)}^2\le2L(f(x_k)-f^\star)$. Minimizing the strong convexity inequality at $x_k$ over the whole space gives $2\mu(f(x_k)-f^\star)\le\norm{\nabla f(x_k)}^2$. Inserting this lower bound into the descent bound gives the contraction factor $1-\kappa^{-1}$.
\end{proof}

The orthogonality $\nabla f(x_{k+1})\perp S_k$ plays the role of the coupling step of Nesterov's method, and every lower bound must respect it. We call a method, or a family of methods indexed by $\kappa$, \emph{uniformly accelerated} if there are absolute constants $C,c>0$ with
\begin{equation}\label{eq:accel}
f(x_k)-f^\star\le C\exp\bigl(-ck/\sqrt\kappa\bigr)(f(x_0)-f^\star)
\end{equation}
for every $\kappa$, every $f\in\Fc_{\mu,L}$ with $L/\mu=\kappa$, every starting point and every $k$. A failure of \eqref{eq:accel} is a worst-case statement. It is compatible with the convex rate $O(1/k^2)$.

SESOP attains this rate and is conjugate gradients on quadratics, and the methods of Section~\ref{sec:filters} that keep its directions inherit both known properties.

\begin{proposition}[orth-method and conjugate gradients \cite{nemirovski1982orth,narkiss2005sequential,hestenes1952methods}]\label{prop:sesop}
Let the iterates satisfy \eqref{eq:step}.
\begin{enumerate}
\item If $S_k$ contains the anchor $x_k-x_0$ and the aggregate $\sigma_k=\sum_{i=0}^kw_i\nabla f(x_i)$ with Nemirovski's weights and $w_0=1$, if $f$ is convex with $L$-Lipschitz gradient and a minimizer $x^\star$, and if the subspace minimizers exist, then $f(x_k)-f^\star\le2L\norm{x_0-x^\star}^2/(k+1)^2$ for $k\ge1$.
\item If $f(x)=\frac12(x-x^\star)^\top A(x-x^\star)$ with $A\succ0$ and $\spn\{\nabla f(x_k),s_k\}\subseteq S_k\subseteq\spn\{\nabla f(x_0),\dots,\nabla f(x_k),s_1,\dots,s_k\}$ for every $k$, then the iterates coincide with those of conjugate gradients started at $x_0$, and $f(x_k)-f^\star\le4\bigl(\frac{\sqrt\kappa-1}{\sqrt\kappa+1}\bigr)^{2k}(f(x_0)-f^\star)$.
\end{enumerate}
\end{proposition}
\begin{proof}
Part~1 is Nemirovski's convex rate for the orth-method \cite{nemirovski1982orth}. Its proof for SESOP \cite[Theorem~1]{narkiss2005sequential} uses three facts: exact minimization makes each gradient $\nabla f(x_i)$ orthogonal to the anchor and to the previous aggregate, each step decreases $f$ at least as much as the gradient step of length $1/L$, and Nemirovski's weights satisfy $w_k^2=\sum_{i\le k}w_i\ge(k+2)^2/4$. With convexity they give $w_k^2(f(x_{k+1})-f^\star)\le\frac L2\norm{x_0-x^\star}^2$, which is the stated bound.

For part~2, the $k$th point of conjugate gradients minimizes $f$ over $x_0$ plus the Krylov space of order $k$, and its step lies in $\spn\{\nabla f(x_k),s_k\}$. If the iterates agree with conjugate gradients up to $k$, then $x_k+S_k$ lies in $x_0$ plus the next Krylov space and contains the next point of conjugate gradients, which is therefore its minimizer. The rate is the standard bound for conjugate gradients \cite[Chapter~5]{nocedal2006numerical}.
\end{proof}

The first tool is a quadratic potential certificate inspired by Nesterov's estimate sequences \cite{nesterov2018lectures}. It shows that an exact step is accelerated whenever a virtual momentum $\pi_k$, which is never computed, lies in the search space, and our global accelerated bounds rest on it. Fix an origin $b$, a number $0<\lambda\le\mu$ that plays the role of $\mu$, $L'\ge L$ and numbers $a_k\in(0,1]$ with $a_k^2/\lambda\le1/L'$ for $k\ge b$ that play the role of $\sqrt{\mu/L}$. With $\rho_k=1-a_k$ and some $\pi_b\in\R^n$ define
\begin{equation}\label{eq:prec}
\pi_{k+1}=\rho_k\bigl(\pi_k+s_{k+1}\bigr)+\frac{a_k}\lambda \nabla f(x_{k+1})\qquad(k\ge b).
\end{equation}
Inexact steps are measured by the defect
\[
\begin{aligned}
E_k&=\Bigl[\ip{\nabla f(x_{k+1})}{\pi_k+s_{k+1}}-\frac\lambda2\norm{\pi_k+s_{k+1}}^2\Bigr]_+,\\
e_b&=0,\qquad e_{k+1}=\rho_k(e_k+a_kE_k).
\end{aligned}
\]

\begin{lemma}[potential certificate]\label{lem:cert}
Let $f\in\Fc_{\mu,L}$ and let $f(x_{k+1})\le f(x_k)-\norm{\nabla f(x_k)}^2/(2L')$ for $k\ge b$. Then for $j\ge0$
\[
f(x_{b+j+1})-f^\star\le\Bigl(\prod_{i=b}^{b+j-1}\rho_i\Bigr)C_b+e_{b+j},\qquad C_b=f(x_b)-f^\star+\frac\lambda2\norm{x_b-\pi_b-x^\star}^2 .
\]
If $\pi_k+s_{k+1}$ lies in a subspace $S$, then $E_k\le\norm{P_S\nabla f(x_{k+1})}^2/(2\lambda)$, where $P_S$ is the orthogonal projector onto $S$. For iterates of \eqref{eq:step} with $L'=L$ and $\pi_k\in S_k$ for $b\le k<b+j$ the defect $e_{b+j}$ vanishes. For $\pi_b=0$, $C_b\le2(f(x_b)-f^\star)$.
\end{lemma}
\begin{proof}
Put $\Delta_k=f(x_k)-f^\star$, $v_k=x_k-\pi_k$ and
\[
V_k=\Delta_{k+1}+\frac\lambda2\norm{v_k-x^\star}^2.
\]
By \eqref{eq:prec}, $v_{k+1}=\rho_kv_k+a_kx_{k+1}-a_k\nabla f(x_{k+1})/\lambda$. Expanding the square and using strong convexity at $x_{k+1}$, $\lambda\le\mu$ and the definition of $E_k$ gives
\[
\begin{aligned}
\frac\lambda2\norm{v_{k+1}-x^\star}^2
&\le\rho_k\frac\lambda2\norm{v_k-x^\star}^2-a_k\Delta_{k+1}\\
&\quad+\frac{a_k^2}{2\lambda}\norm{\nabla f(x_{k+1})}^2+a_k\rho_kE_k.
\end{aligned}
\]
Adding the descent bound for $\Delta_{k+2}$ and using $a_k^2/\lambda\le1/L'$ yields
\[
V_{k+1}\le\rho_kV_k+a_k\rho_kE_k.
\]
Since $V_b\le C_b$, iteration gives $V_{b+j}\le(\prod_{i=b}^{b+j-1}\rho_i)C_b+e_{b+j}$, and $\Delta_{b+j+1}\le V_{b+j}$ proves the first claim.

If $\pi_k+s_{k+1}\in S$, the bracket in $E_k$ is at most $\max_{u\in S}\{\ip{P_S\nabla f(x_{k+1})}u-\frac\lambda2\norm u^2\}=\norm{P_S\nabla f(x_{k+1})}^2/(2\lambda)$. For iterates of \eqref{eq:step} Lemma~\ref{lem:descent} gives the descent assumption with $L'=L$, $s_{k+1}\in S_k$ and $\nabla f(x_{k+1})\perp S_k$, so $\pi_k\in S_k$ gives $\ip{\nabla f(x_{k+1})}{\pi_k+s_{k+1}}=0$ and $E_k=0$. For $\pi_b=0$, $\frac\lambda2\norm{x_b-x^\star}^2\le\frac\mu2\norm{x_b-x^\star}^2\le f(x_b)-f^\star$.
\end{proof}

The second tool builds the hard functions. The lower bounds prescribe points, gradients and values and extend them to a function of the class. For $0<\mu<L$ and data $(x_i,g_i,f_i)$ the ordered interpolation residual is
\begin{equation}\label{eq:residual}
C_{ij}=f_i-f_j-\ip{g_j}{x_i-x_j}-\frac{\tfrac1L\norm{g_i-g_j}^2+\mu\norm{x_i-x_j}^2-\tfrac{2\mu}{L}\ip{g_i-g_j}{x_i-x_j}}{2(1-\mu/L)}.
\end{equation}
For finite data, $C_{ij}\ge0$ is the interpolation condition of \cite[Theorem~4]{taylor2017smooth}. The extension to countable bounded data is routine, and we include its proof for completeness.

\begin{lemma}[bounded countable interpolation]\label{lem:interp}
Let $0<\mu<L$. Consider a nonempty finite or countable bounded set of data
\[
\{(x_i,g_i,f_i)\}_{i\in I}\subseteq\R^n\times\R^n\times\R
\]
with $C_{ij}\ge0$ for all $i,j\in I$. Then there is $F\in\Fc_{\mu,L}(\R^n)$ with $F(x_i)=f_i$ and $\nabla F(x_i)=g_i$ for every $i\in I$. If the triple $(0,0,0)$ belongs to the data, then $0$ is the minimizer of $F$ and $F^\star=0$.
\end{lemma}

\begin{proof}
As in \cite{taylor2017smooth}, put $a=L-\mu>0$, $q_i=g_i-\mu x_i$ and $u_i=\ip{q_i}{x_i}-f_i+\frac\mu2\norm{x_i}^2$. The conditions $C_{ij}\ge0$ are equivalent to
\begin{equation}\label{eq:shifted}
u_j\ge u_i+\ip{x_i}{q_j-q_i}+\frac1{2a}\norm{q_j-q_i}^2\qquad(i,j\in I).
\end{equation}
Define
\[
\Psi(y)=\sup_{i\in I}\Bigl\{u_i+\ip{x_i}{y-q_i}+\frac1{2a}\norm{y-q_i}^2\Bigr\}.
\]
Boundedness of the data makes $\Psi$ finite and continuous, since subtracting $\norm y^2/(2a)$ leaves a supremum of affine functions with uniformly bounded slopes and intercepts. Thus $\Psi$ is $1/a$-strongly convex. By \eqref{eq:shifted}, $\Psi(q_j)=u_j$ and the $j$th quadratic supports $\Psi$ there, so $x_j\in\partial\Psi(q_j)$. Conjugacy \cite[Proposition~12.60]{rockafellar1998variational} gives $\nabla\Psi^*(x_j)=q_j$, $\Psi^*(x_j)=\ip{x_j}{q_j}-u_j$, and an $a$-Lipschitz gradient for $\Psi^*$. Hence $F=\Psi^*+\frac\mu2\norm{\cdot}^2$ belongs to $\Fc_{\mu,L}$ and interpolates all values and gradients. If $(0,0,0)$ is among the data, then $\nabla F(0)=F(0)=0$.
\end{proof}

\section{When memory blocks acceleration}\label{sec:main}

By Lemma~\ref{lem:cert} a method is accelerated once its subspaces contain the virtual momentum $\pi_k$, an average over roughly the last $\sqrt\kappa$ iterations. The lower bounds below constrain the window and the forgetting rate of memories fixed in advance. A memory whose weights do not decay on it keeps directions far older than the horizon (Section~\ref{sec:slow}), and one that forgets much faster sees only its last few gradients (Section~\ref{sec:fastpair}). Section~\ref{sec:lower} combines the two bounds, and Section~\ref{sec:smooth} shows that smoothness does not help.

\subsection{Memories fixed in advance}\label{sec:memories}

SESOP combines a window of length one with the anchor $x_k-x_0$, which sums all past steps, and with weighted sums of all past gradients, so it never forgets anything.

\begin{definition}[windows with aggregates]\label{def:separable}
A memory with an integer window length $m\ge0$ and aggregates consists of the last $m$ gradients and steps, finitely many gradient sums $\sigma^a_k=\sum_{i=0}^{k}w_a(i)\nabla f(x_i)$ and finitely many step sums $\theta^b_k=\sum_{i=0}^{k}v_b(i)s_i$. The weight sequences are fixed before the run and do not depend on the objective. An exact subspace method uses this memory if for every $k\ge0$
\[
\nabla f(x_k)\in S_k\subseteq\spn\bigl(\{\nabla f(x_i),s_i:k-m<i\le k\}\cup\{\sigma^1_k,\dots,\sigma^d_k,\theta^1_k,\dots,\theta^{d'}_k\}\bigr).
\]
The choice of $S_k$ inside this space may depend on the points, values and gradients observed at $x_0,\dots,x_k$.
\end{definition}

Multiplying a register by a nonzero scalar does not change its direction, so forgetting is measured on normalized weights. A linear register $Z_k=Z_{k-1}A_k+[\nabla f(x_k),s_k]B_k$ with both matrix sequences $A_k,B_k$ fixed in advance, each $A_k$ invertible, and $Z_{-1}=0$ satisfies $Z_k\Phi_k^{-1}=\sum_{i=0}^{k}[\nabla f(x_i),s_i]B_i\Phi_i^{-1}$ with $\Phi_k=A_0\cdots A_k$, so its span lies in the span of gradient and step sums whose normalized weights are the entries of $B_i\Phi_i^{-1}$. The raw weights of the geometric filter $z_k=\rho z_{k-1}+\nabla f(x_k)$ decay, while its normalized weights $\rho^{-i}$ grow.

The \emph{growth rate} $\Gamma$ and the \emph{forgetting index} $\delta$ of the memory are
\[
\Gamma=\max\Bigl\{1,\ \max_a\limsup_{i\to\infty}|w_a(i)|^{1/i},\ \max_b\limsup_{i\to\infty}|v_b(i)|^{1/i}\Bigr\},\qquad\delta=1-\frac1\Gamma .
\]
A maximum over an empty family is omitted, so a pure window has $\Gamma=1$ and $\delta=0$. A geometric filter with factor $\rho\in(0,1]$ has $\delta=1-\rho$. We call the memory \emph{slowly forgetting} if $\Gamma=1$, that is, if all its weights grow subexponentially.
GD\nobreakdash-ELS, the plane rule, windows and SESOP with bounded, polynomial or Nemirovski's weights are slowly forgetting. A filter with factor $1-1/\kappa$ forgets exponentially, but $\delta=1/\kappa$ is small on the scale of acceleration. 

\subsection{Remembering too much}\label{sec:slow}

Our first lower bound concerns memories with a window of length $m\le p-2$ and weights that grow more slowly than $1/r$, where $r=1-p/\kappa$. For such a memory one function makes every method that uses it follow the same trajectory, on which the gap shrinks only by the factor $r^2$ per step.

\begin{theorem}[slowly forgetting memories]\label{thm:slow}
There is an absolute constant $C_0\ge1$ with the following property. Fix a memory of Definition~\ref{def:separable} with window $m$ and forgetting index $\delta$. Let $p\ge\max(8,m+2)$ be an integer, $\kappa\ge C_0p^2$ and $\delta<p/\kappa$. Then there are a dimension $n$, a function $F\in\Fc_{1/\kappa,1}(\R^n)$ and a starting point such that every exact subspace method \eqref{eq:step} using this memory generates the same trajectory on $F$, and $F(x_{k+1})-F^\star=(1-p/\kappa)^2(F(x_k)-F^\star)$ for all large $k$. Hence $\liminf_{\varepsilon\downarrow0}N(\varepsilon)/\log(1/\varepsilon)\ge(\kappa-p)/(2p)\ge3\kappa/(8p)$. The same holds for $m=1$, $p=3$ and every $\kappa\ge12$.
\end{theorem}

The proof of Theorem~\ref{thm:slow} is given at the end of this subsection.

\begin{corollary}[SESOP and conjugate gradients]\label{cor:sesop}
Let SESOP use any finite number of aggregates with subexponential weights, in particular Nemirovski's or the linear weights, with or without the previous step and the anchor. For every $\kappa\ge12$ there is $F\in\Fc_{1/\kappa,1}$ on which $\liminf_{\varepsilon\downarrow0}N(\varepsilon)/\log(1/\varepsilon)\ge(\kappa-3)/6\ge\kappa/8$. The plane form of conjugate gradients and GD\nobreakdash-ELS follow the same trajectory. None of these methods is uniformly accelerated.
\end{corollary}
\begin{proof}
These memories have the window $m=1$ and subexponential weights: the anchor is the step sum with weights $1$, and Nemirovski's weights satisfy $w(k-1)+\frac12\le w(k)\le w(k-1)+1$. Hence $\Gamma=1$, and Theorem~\ref{thm:slow} with $m=1$ and $p=3$ applies to every rule using this memory, which includes all methods of the corollary. By \eqref{eq:accel} a uniformly accelerated method has $\limsup_{\varepsilon\downarrow0}N(\varepsilon)/\log(1/\varepsilon)\le\sqrt\kappa/c$, which is below $\kappa/8$ for large $\kappa$.
\end{proof}

The trajectory contracts by $(1-3/\kappa)^2$ per step, so the exact asymptotic coefficient is $[-2\log(1-3/\kappa)]^{-1}\sim\kappa/6$ as $\kappa\to\infty$.

For the proof we prescribe points $x_k$, gradients $g_k$ (which become $\nabla F(x_k)$) and values $f_k$ along an infinite trajectory together with the minimizer $0$, and Lemma~\ref{lem:interp} extends the data to a function of the class if all residuals are nonnegative. If every step is a negative multiple of the current gradient and every next gradient is orthogonal to the whole memory, the method with the largest admissible spaces generates the trajectory, and by the next lemma so does every method using the memory. The trajectory consists of a cyclic tail (Lemma~\ref{lem:cyclic}) and an orthogonal initial segment (Lemma~\ref{lem:prefix}), which Lemma~\ref{lem:bridge} joins, and Lemma~\ref{lem:generic-t} chooses the step lengths of the segment.

\begin{lemma}[transfer]\label{lem:transfer}
Let $f$ be differentiable and strongly convex, and let an exact subspace method with spaces $S_k$ generate a trajectory with $s_{k+1}\in\spn\{\nabla f(x_k)\}$ for all $k$. Let another rule produce, along the same history, spaces with $\spn\{\nabla f(x_k)\}\subseteq\tilde S_k\subseteq S_k$. Then the second rule, started at the same $x_0$, generates the same trajectory.
\end{lemma}

\begin{proof}
Suppose the histories agree up to $k$. Since $s_{k+1}\in\spn\{\nabla f(x_k)\}\subseteq\tilde S_k$, the unique minimizer over $x_k+S_k$ belongs to $x_k+\tilde S_k\subseteq x_k+S_k$. It is therefore also the unique minimizer over the smaller space. Induction proves the claim.
\end{proof}

Throughout this subsection and the next $L=1$, $\bar\mu=1-\mu$, $C_{ij}$ is the residual \eqref{eq:residual}, and the index $\star$ refers to the minimizer $(0,0,0)$. The slow part of the trajectory is a cyclic tail. Its gradient directions run through $p$ orthonormal vectors, so that gradients at the lags $1,\dots,p-1$ are orthogonal, and their lengths decay by the factor $r$ per step. With the right step length $\tau$ and values $f_k$ all off-diagonal residuals are bounded below by a fixed multiple of $r^{2\min(i,j)}$, a margin that later absorbs perturbations.

Both tails below are orbits of an orthogonal map $Q$ with a factor $r\in(0,1)$,
\begin{equation}\label{eq:orbit}
(x_n,g_n,f_n)=\bigl(r^nQ^nX,\,r^nQ^nG,\,\tfrac c2r^{2n}\bigr)\quad(n\ge0),\qquad\ip GX=c .
\end{equation}
A common shift of indices by $i$ multiplies all their residuals, including those with the minimizer, by $r^{2i}$, and $\ip GX=c$ makes both residuals of $(X,G,\frac c2)$ with the minimizer equal to $C_\star:=C_{0\star}=C_{\star0}$. For $\ell>0$ the ordered interpolation residuals defined in \eqref{eq:residual} between $(X,G,\frac c2)$ and $(tQ^\ell X,tQ^\ell G,\frac c2t^2)$ are polynomials of degree at most two in $t$, equal $C_\star$ at $t=0$ and, again by $\ip GX=c$, have the coefficient $C_\star$ at $t^2$. With $t=r^\ell$ this gives
\begin{equation}\label{eq:orbit-res}
C_{0\ell}=C_\star(1-t)^2+tK_\ell^+,\qquad C_{\ell0}=C_\star(1-t)^2+tK_\ell^-,
\end{equation}
where $K_\ell^+$ and $K_\ell^-$ are the residuals between $(X,G,\frac c2)$ and $(Q^\ell X,Q^\ell G,\frac c2)$ in the same order.

\begin{lemma}[cyclic tails]\label{lem:cyclic}
Let either (i) $p=3$, $0<\mu\le\frac1{12}$, $r=1-3\mu$ and $c=[3\mu(1-2\mu)]^{-1}$, or (ii) $p\ge8$, $\kappa\ge C_0p^2$ for a sufficiently large absolute constant $C_0$, $\mu=1/\kappa$, $r=1-p/\kappa$ and $c=\kappa/p+3p/4$. Put $\tau=c(1-r^p)$, $\zeta_n=e_{1+(n\bmod p)}\in\R^p$ and
\[
g_n=r^n\zeta_n,\qquad x_n=\tau\sum_{j=n}^{\infty}g_j,\qquad f_n=\frac c2r^{2n}\qquad(n\ge0).
\]
Then $x_{n+1}-x_n=-\tau g_n$, $C_{ij}\ge\eta r^{2\min(i,j)}$ for $i\ne j$ and $C_{i\star},C_{\star i}\ge\eta r^{2i}$ for some $\eta>0$.
\end{lemma}

\begin{proof}
All residuals reduce to the scalars $\mathcal A$ and $\mathcal E_q$ defined below, and the lemma follows once they are positive.

\emph{Step 1: reduction.} Let $R$ be the cyclic shift of $\R^p$. Grouping the series by residues gives $x_n=r^nR^nx_0$ with $x_0=c(1,r,\dots,r^{p-1})^\top$, and the first coordinate of $rRx_0-x_0$ is $c(r^p-1)$ while the others vanish, which gives the step identity. Put
\[
S=\sum_{i=0}^{p-1}r^{2i},\qquad T_d=\sum_{i=0}^{p-1}r^ir^{(i-d)\bmod p},\qquad\mathcal A=(1+\mu)c-1-\mu c^2S
\]
and, for $1\le q\le\lfloor p/2\rfloor$,
\[
\mathcal E_q=c\bigl[(1-r^q)+\mu(1-r^{p-q})\bigr]-1-\mu c^2(S-T_q).
\]
The data form an orbit \eqref{eq:orbit} with $Q=R$, $X=x_0$ and $G=e_1$, and $C_\star=\mathcal A/(2\bar\mu)$. Fix a lag $j>0$, $t=r^j$ and $d=j\bmod p$. If $d=0$, then $K_j^\pm=0$. If $d\ne0$ and $q=\min(d,p-d)$, the identities $\norm{e_1-R^de_1}^2=2$, $\norm{x_0-R^dx_0}^2=2c^2(S-T_d)$ and $\ip{e_1-R^de_1}{x_0-R^dx_0}=c(2-r^d-r^{p-d})$ give $\min K_j^\pm=\mathcal E_q/\bar\mu$. Once $\mathcal A$ and all $\mathcal E_q$ are positive, \eqref{eq:orbit-res} has a positive minimum on $0\le t\le r^p$ for $d=0$ and on $0\le t\le1$ for $d\ne0$. The finitely many residue classes therefore give a common margin $\eta>0$.

\emph{Step 2: case (i).} Inserting $r$ and $c$ gives
\[
\mathcal A=\frac{(1-3\mu)(9\mu^2-5\mu+2)}{3(1-2\mu)^2}>0,\qquad\mathcal E_1=\frac{\mu(1+\mu)(1-3\mu)}{(1-2\mu)^2}>0 .
\]

\emph{Step 3: case (ii), uniformly in $p$.} Put $u=p^2\mu$, so $r=1-u/p$ and $c=p(1/u+3/4)$. Taylor's formula for $r^j$, $j\le2p$, gives $S=p-(p-1)u+O(pu^2)$. Using
\[
S-T_1=\frac{(1-r^{p-1})(1-r^p)}{1+r},\qquad
S-T_q=\frac12\sum_{i=0}^{p-1}(r^i-r^{(i-q)\bmod p})^2\le\frac{u^2q(p-q)}{2p},
\]
where the inequality follows from $|r^i-r^j|\le|i-j|u/p$, we obtain
\[
\begin{aligned}
\mathcal A&=p/4-2+1/p+O(pu),\\
\mathcal E_1&=u\bigl(1/4+3/(2p)-1/p^2+O(u)\bigr),\\
\mathcal E_q&=q-1+O(uq)\qquad(2\le q\le p/2).
\end{aligned}
\]
All remainder constants are absolute, since $j\le2p$ bounds the normalized Taylor remainders uniformly for $0<u\le1/2$. The leading terms of $\mathcal A/p$, $\mathcal E_1/u$ and $\mathcal E_q/q$ are at least $1/64$, $1/4$ and $1/2$, respectively. Thus all are positive for $u\le1/C_0$ with an absolute $C_0$.
\end{proof}

Steps are multiples of gradients, so the window at iteration $k$ is spanned by $g_{k-m},\dots,g_k$, whose lags to $k+1$ lie between $1$ and $m+1\le p-1$. The next tail gradient is therefore orthogonal to the window. Thus Lemmas~\ref{lem:interp}, \ref{lem:transfer} and~\ref{lem:cyclic} already give the required lower bound when there are no aggregates or anchor.

To include cumulative sums, we need two further lemmas. A gradient sum $\sum_{i\le k}w(i)g_i$ is generally not orthogonal to the next tail gradient. If the weights grow more slowly than $1/r$, its tail is geometrically small, and a small perturbation of the tail directions makes every tail gradient orthogonal to the sums if they start from suitable values. The next lemma constructs an initial segment of orthogonal gradient steps with a positive interpolation margin; the following lemma joins it to the perturbed tail.

\begin{lemma}[orthogonal initial segment]\label{lem:prefix}
For some absolute constants $c_{\mathrm{pre}}>0$ and $\epsilon_0\in(0,1/10]$, the following holds. Let $0<\mu\le1/12$, $B\ge1$, and let $G_0,\dots,G_{B-1}$ be orthogonal with $\norm{G_i}^2=10^{-i}$. For $t_i\in[2-\epsilon_0,2+\epsilon_0]$, define $X_i=\sum_{h=i}^{B-1}t_hG_h$ and $F_i=\sum_{h=i}^{B-1}\norm{G_h}^2$ for $i<B$, and $(X_B,G_B,F_B)=(0,0,0)$. Then $C_{ij},C_{ji}\ge c_{\mathrm{pre}}\norm{G_i}^2$ for $i<j\le B$.
\end{lemma}

\begin{proof}
At $t_i=2$, orthogonality and geometric summation show that the nonquadratic parts of $C_{ij}$ and $C_{ji}$, divided by $\norm{G_i}^2$, are at least $1$ and $8/9$. Their common quadratic correction, normalized in the same way, is at most
\[
\frac{1+1/10+4\mu/9}{2(1-\mu)}\le\frac23.
\]
Hence both residuals are at least $(2/9)\norm{G_i}^2$. Perturbing all $t_h$ by at most $\epsilon$ changes each normalized residual by $O(\epsilon)$, uniformly in $B,i,j$ and $\mu\le1/12$, because the normalized coefficients are controlled by $\sum_{h\ge0}10^{-h}=10/9$. A sufficiently small absolute $\epsilon_0$ therefore preserves a positive margin $c_{\mathrm{pre}}$.
\end{proof}

The next lemma joins such a segment to a perturbed tail, and its proof is the technical core of this section. The unknowns are a small perturbation $h_n$ of every tail direction and finitely many coordinates of the initial gradients, and every tail gradient has to be orthogonal to the anchor and to each aggregate. These are infinitely many nonlinear conditions, but at the unperturbed point they decouple. The lag conditions form a cyclic linear system, which can be inverted explicitly, and the conditions for the anchor and the aggregates fix the components of the perturbation along the initial segment through an invertible matrix $\mathbf A$. The implicit function theorem then applies in a space of geometrically weighted sequences, and this is where the growth of the weights enters. It needs a small parameter $\vartheta$, and since the initial segment has size of order $\vartheta^{-1}$, the slow tail starts only after the gap has decreased by a large factor, which we did not try to control.

\begin{lemma}[cyclic bridge]\label{lem:bridge}
Let the data of Lemma~\ref{lem:cyclic} be given and $1\le m\le p-2$. Let $B\ge1$, $t\in(2-\epsilon_0,2+\epsilon_0)^B$, with $\epsilon_0$ from Lemma~\ref{lem:prefix}, and let $w_1,\dots,w_{B-1}$ be real sequences with $\omega=\max_a\limsup_{i\to\infty}|w_a(i)|^{1/i}<1/r$ (with $\omega=0$ if $B=1$) such that the $B\times B$ matrix $\mathbf A$ with the rows $t,w_1|_{[0,B)},\dots,w_{B-1}|_{[0,B)}$ is invertible. Then there are $F\in\Fc_{\mu,1}(\R^{p+B})$ and a trajectory of the exact method with the spaces $\spn\{g_i,s_i:k-m<i\le k\}+\spn\{x_k-x_0,\sigma^1_k,\dots,\sigma^{B-1}_k\}$, $\sigma^a_k=\sum_{i=0}^kw_a(i)g_i$, whose steps are $-t_kg_k$ for $k<B$ and $-\tau g_k$ for $k\ge B$, and $F(x_k)-F^\star=\frac c2r^{2k}$ for $k\ge B$.
\end{lemma}

\begin{proof}
\emph{Joint system.} Choose $\gamma$ with $\max\{1,\omega\}<\gamma<1/r$ and $\beta\in(r\gamma,1)$, so that $|w_a(j)|\le C_a\gamma^j$ for all $a$ after enlarging the constants. The tail sums $\sum_{j=n}^{\infty}w_a(j)g_j$ are then $O((r\gamma)^n)=O(\beta^n)$. In the space $E=\R^B$ orthogonal to $\R^p$ take the vectors $G_i$ of Lemma~\ref{lem:prefix}. For triangular parameters $\lambda_{ij}$, $j<i<B$, put
\[
\begin{aligned}
Z_i(\lambda)&=G_i+\sum_{j=0}^{i-1}\lambda_{ij}G_j,\qquad U(\lambda)=\sum_{i=0}^{B-1}t_iZ_i(\lambda),\\
V_a(\lambda)&=\sum_{i=0}^{B-1}w_a(i)Z_i(\lambda)\qquad(1\le a<B).
\end{aligned}
\]
For $n\ge B$ seek $v_n=\zeta_n+h_n$ and set $g_n=r^nv_n$, $x_n=\tau\sum_{j=n}^{\infty}g_j$ and $f_n=\frac c2r^{2n}$. For every $n\ge B$ we impose the $p+B$ equations
\begin{align*}
&\norm{v_n}=1,\quad\ip{v_n}{v_{n+t}}=0\ \ (1\le t\le p-1)&&\text{(norm and lag conditions)},\\
&\ip{v_n}{U(\lambda)}=\vartheta\ip{v_n}{x_n}&&\text{(anchor condition)},\\
&\textstyle\ip{v_n}{V_a(\lambda)}=\vartheta\bigl\langle v_n,\sum_{j\ge n}w_a(j)g_j\bigr\rangle\ \ (1\le a<B)&&\text{(aggregate conditions)}.
\end{align*}
Define correction vectors by $\sum_i\mathbf A_{\ell i}T_i=c_\ell$ for $0\le\ell<B$, where $c_0=x_B$ and $c_a=\sum_{n\ge B}w_a(n)g_n$ for $a\ge1$, and put $P_i=Z_i(\lambda)-\vartheta T_i$. Finally impose the finite Gram equations $\ip{P_i}{P_j}=0$ for $0\le j<i<B$.

\emph{Function spaces.} The unknowns are $\lambda\in\R^{B(B-1)/2}$ and $h=(h_n)_{n\ge B}$ in the space $\ell^\infty_\beta$ of sequences in $\R^{p+B}$ with the norm $\norm h_\beta=\sup_n\beta^{-n}\norm{h_n}$. Since $\zeta_n\perp E$ and $r<r\gamma<\beta$, the residuals of the tail equations at index $n$ are $O(\beta^n)$ for $h\in\ell^\infty_\beta$, so they define a map into $\ell^\infty_\beta$, and the Gram equations a map into $\R^{B(B-1)/2}$. Since $|w_a(n)|r^n\le C(\gamma r)^n$ and $\gamma r<\beta$, the vector $x_B$ and the sums $\sum_{n\ge B}w_a(n)g_n$ depend affinely and boundedly on $h$, and so do the corrections $T_i$. The whole map is therefore continuously differentiable near $(0,0,0)$.

\emph{Invertibility.} The system holds at $(h,\lambda,\vartheta)=(0,0,0)$. At this point the derivative with respect to $(h,\lambda)$ is block diagonal. The Gram equations do not depend on $h$ at $\vartheta=0$, because then $P_i=Z_i(\lambda)$. The tail equations have zero derivative with respect to $\lambda$ at $h=0$, because $U(\lambda)$ and $V_a(\lambda)$ lie in $E$, which is orthogonal to every $\zeta_n$. In the cyclic block the linearized lag equations read $a^t_n+a^{p-t}_{n+t}=q^t_n$ with $a^t_n=\ip{h_n}{\zeta_{n+t}}$ and right-hand sides $q^t_n$, and their inverse is
\[
a^t_n=\sum_{\ell=0}^{\infty}\bigl(q^t_{n+\ell p}-q^{p-t}_{n+t+\ell p}\bigr),
\]
which telescopes and is unique among decaying sequences. Defining
\[
\norm q_\beta=\max_{1\le t<p}\sup_{n\ge B}\beta^{-n}|q_n^t|
\]
gives $\norm{a^t}_\beta\le\frac{1+\beta^t}{1-\beta^p}\norm q_\beta$.
The linearized norm equation fixes the remaining cyclic coordinate $\ip{h_n}{\zeta_n}$ directly. Thus the inverse of the cyclic block is bounded on $\ell^\infty_\beta$. The block of the extra coordinates consists of the $B$ scalar products with $U(0),V_1(0),\dots,V_{B-1}(0)$. In the orthonormal basis $G_i/\norm{G_i}$ of $E$ its matrix is $\mathbf A\operatorname{diag}(\norm{G_i})$, which is invertible. The finite block is diagonal with entries $\norm{G_j}^2>0$. The implicit function theorem gives a solution with $\norm{h}_\beta+\norm{\lambda}=O(\vartheta)$.

\emph{Actual initial segment.} For $i<B$ put $g_i=\vartheta^{-1}P_i$, $x_i=x_B+\vartheta^{-1}\sum_{h=i}^{B-1}t_hP_h$ and $f_i=\vartheta^{-2}F_i$. The corrections give $x_0=\vartheta^{-1}U(\lambda)$ and $\sum_{i=0}^{\infty}w_a(i)g_i=\vartheta^{-1}V_a(\lambda)$ for every row $1\le a<B$. The initial gradients are mutually orthogonal and each initial step is $-t_ig_i$. The tail equations then give $\ip{g_n}{x_n-x_0}=0$ and $\ip{g_n}{\sum_{i=0}^{n-1}w_a(i)g_i}=0$ for $n\ge B$ and all $1\le a<B$. For $B\le n\le B+p-2$ the lag conditions make $g_n$ orthogonal to the tail gradients with indices from $B$ to $n-1$, which leaves the $B$ relations $\ip{g_n}{\sum_{i=0}^{B-1}\mathbf A_{\ell i}g_i}=0$, one for each row $\ell$ of $\mathbf A$. Since $\mathbf A$ is invertible, $g_n\perp g_i$ for all $i<B$. This range suffices, because for $m\le p-2$ the window at iteration $n-1$ contains an initial gradient only if $n\le B+m$.

\emph{Optimality.} For $k+1<B$ every generator of $S_k$ is a combination of $g_0,\dots,g_k$, and $g_{k+1}$ is orthogonal to them. For $k+1\ge B$ the lag conditions and the relations above make $g_{k+1}$ orthogonal to the window. The anchor relation gives $\ip{g_{k+1}}{x_k-x_0}=\ip{g_{k+1}}{x_{k+1}-x_0}-\ip{g_{k+1}}{s_{k+1}}=0$, and the aggregate relations give $g_{k+1}\perp\sigma^a_k$. Thus $s_{k+1}\in S_k$ and $g_{k+1}\perp S_k$, and by strong convexity $x_{k+1}$ is the unique minimizer over $x_k+S_k$.

\emph{Interpolation.} The tail points and gradients are $O(r^n)$, their perturbations are $O(\vartheta(r\beta)^n)$, and the values are unchanged. Since the remaining terms of each residual are quadratic, $|C_{ij}-C_{ij}^0|\le M\vartheta r^{2\min(i,j)}$ uniformly for $i,j\ge B$, with the corresponding bound $M\vartheta r^{2i}$ for pairs with the minimizer. Taking $M\vartheta<\eta/2$ preserves all tail margins from Lemma~\ref{lem:cyclic}. For off-diagonal pairs that involve the initial segment we scale the data by $(x,g,f)\mapsto(\vartheta x,\vartheta g,\vartheta^2f)$, which multiplies all residuals by $\vartheta^2$. The scaled initial data converge to $(X_i,G_i,F_i)$ of Lemma~\ref{lem:prefix} and the scaled tail data uniformly to $(0,0,0)=(X_B,G_B,F_B)$, so the bound $c_{\mathrm{pre}}\norm{G_i}^2$ of Lemma~\ref{lem:prefix} and uniform continuity of the residual keep these residuals positive for small $\vartheta$. Lemma~\ref{lem:interp} completes the construction.
\end{proof}

Lemma~\ref{lem:bridge} handles the anchor and gradient sums. Along a trajectory with steps $s_{i+1}=-t_ig_i$ a step sum $\theta^b_k=\sum_{i\le k}v_b(i)s_i$ equals $-\sum_{j=0}^{k-1}t_j\tilde v_b(j)g_j$ with $\tilde v_b=v_b(\cdot+1)$, a gradient sum without its last term, and the anchor is the case $v_b\equiv1$. Its initial weights depend on the lengths $t_j$, chosen so that the initial segment fixes these sums.

\begin{lemma}[generic step lengths]\label{lem:generic-t}
Let $w_1,\dots,w_e$, $\tilde v_0\equiv1,\tilde v_1,\dots,\tilde v_f$ be sequences and $\tau>0$. There are $B\ge1$ and $t_0,\dots,t_{B-1}\in I=(2-\epsilon_0,2+\epsilon_0)$, with $\epsilon_0$ from Lemma~\ref{lem:prefix}, such that, with $t_i=\tau$ for $i\ge B$ and $(t\tilde v_c)(i)=t_i\tilde v_c(i)$, the restriction to $[0,B)$ is injective on $\mathcal U=\spn\{w_1,\dots,w_e,t\tilde v_0,\dots,t\tilde v_f\}$.
\end{lemma}

\begin{proof}
For a positive sequence $t$ the dimension of $\mathcal U$ is at most $e+f+1$. Let $t^\ast$ attain its maximum $r$, choose $r$ linearly independent members of the family for $t^\ast$ and indices $i_1,\dots,i_r$ at which their values form a nonsingular matrix. Its determinant is a polynomial in $t_{i_1},\dots,t_{i_r}$ that does not vanish at $t^\ast$, so it does not vanish at some point of $I^r$. Take $B>\max_ji_j$, these values at $i_1,\dots,i_r$, arbitrary values in $I$ at the other indices below $B$, and $t_i=\tau$ for $i\ge B$. Then $\dim\mathcal U\le r$ by maximality, while the restrictions to $[0,B)$ have rank at least $r$, so the restriction is injective.
\end{proof}

\begin{proof}[of Theorem~\ref{thm:slow}]
If $m=0$, enlarge the window to length one; every method using the original memory also uses this enlarged memory, so a lower bound for the latter suffices. We may therefore assume $m\ge1$.
Without aggregates Lemma~\ref{lem:interp} applied to the data of Lemma~\ref{lem:cyclic} gives the instance, as explained after Lemma~\ref{lem:cyclic}. With aggregates we proceed in three steps, using the data of Lemma~\ref{lem:cyclic}, case (i) for $m=1$, $p=3$ and case (ii) otherwise.

\emph{The weights.} Apply Lemma~\ref{lem:generic-t} to the weights $w_a$, to $\tilde v_b$ and to $\tilde v_0\equiv1$, and choose a basis $t,u_1,\dots,u_q$ of $\mathcal U$. Its restrictions to $[0,B)$ are linearly independent, and completing them by rows supported on $[0,B)$ gives an invertible matrix $\mathbf A$ as in Lemma~\ref{lem:bridge}. For $i\ge B$ every $u_j(i)$ is a fixed combination of $w_a(i)$, $\tau\tilde v_b(i)$ and $\tau$, so $\delta<p/\kappa$, that is $\Gamma<1/(1-p/\kappa)$, gives the growth condition of Lemma~\ref{lem:bridge}.

\emph{The trajectory.} Apply Lemma~\ref{lem:bridge} to these data. On its trajectory $g_{k+1}$ is orthogonal to the window, to the anchor and to the sums with the weights $u_j$. The weights $w_a$ and $t\tilde v_b$ are combinations of $t,u_1,\dots,u_q$ as sequences, so $g_{k+1}$ is orthogonal to every gradient sum and, since $g_{k+1}\perp g_k$, to every step sum. Every step is a multiple of the current gradient, so by Lemma~\ref{lem:transfer} every rule using the memory generates this trajectory, and $F(x_k)-F^\star=\frac c2r^{2k}$ for $k\ge B$.

\emph{The rate.} Fix the constructed function $F$ and put $\Delta_0=F(x_0)-F^\star>0$. Since $F(x_k)-F^\star=\frac c2r^{2k}$ for $k\ge B$ (with $B=0$ for a pure window), monotonicity gives, for all sufficiently small $\varepsilon>0$,
\[
N(\varepsilon)=\left\lceil\frac{\log(c/(2\varepsilon\Delta_0))}{-2\log r}\right\rceil.
\]
Thus the fixed initial segment changes only the additive constant in $N(\varepsilon)$. In both cases $p/\kappa\le\frac14$, and $-\log(1-u)\le u/(1-u)$ for $0<u<1$ yields
\[
\lim_{\varepsilon\downarrow0}\frac{N(\varepsilon)}{\log(1/\varepsilon)}
=\frac1{-2\log r}\ge\frac{\kappa-p}{2p}\ge\frac{3\kappa}{8p}.
\tag*{\qed}
\]
\let\qed\relax
\end{proof}

\subsection{Forgetting too fast}\label{sec:fastpair}

A pair of filters with a fixed factor $\rho$ has the normalized weights $\rho^{-i}$, and the construction of Section~\ref{sec:slow} covers it only if $1-\rho<p/\kappa$. Forgetting too fast is just as bad, because such a pair is no better than a short window: with $\theta=1-\rho^{m+2}$, some function costs every method using the window and the pair about $\kappa\theta/m$ iterations per unit of $\log(1/\varepsilon)$. The perturbation of Section~\ref{sec:slow} fails here, because the normalized weights grow faster than the tail decays. We therefore choose the tail with exact orthogonality to the pair, and only \eqref{eq:orbit-res} and Lemmas~\ref{lem:interp} and~\ref{lem:transfer} carry over.

\begin{theorem}[a fast pair]\label{thm:fast}
There is an absolute constant $C_1>0$ with the following property. Let the memory consist of the window of length $m\ge1$ and one pair $z_{k+1}=\rho z_k+\nabla f(x_{k+1})$, $y_{k+1}=\rho(y_k+s_{k+1})$, $z_0=y_0=0$, with $0<\rho<1$. Let $d=m+2$ and $\theta=1-\rho^d$, and assume $\kappa\theta^2\ge C_1d^4$. Then the conclusion of Theorem~\ref{thm:slow} holds with the factor $(1-8d/(\kappa\theta))^2$, and $\liminf_{\varepsilon\downarrow0}N(\varepsilon)/\log(1/\varepsilon)\ge\kappa\theta/(32d)$.
\end{theorem}

\begin{proof}
The tail directions combine two orthonormal families of $d$ vectors, with one component changing sign after $d$ steps. In the product of the next gradient with the filter only the lags that are multiples of $d$ survive, and they come with alternating signs. The weight $\alpha$ of the two components is chosen so that the even and the odd multiples of $d$ cancel exactly, see \eqref{eq:ef-orth} below. The window sees only the lags $1,\dots,d-1$, at which the tail directions are orthogonal anyway, so $d=m+2$ is the smallest admissible block length.

\emph{Parameters.} Put $\mu=1/\kappa$, $\eta=8d/(\kappa\theta)$, $r=1-\eta$, $u=r^d$, $\nu=(\rho/r)^d$,
\[
\alpha=\frac{1-\nu}2,\qquad\tau=\frac{4d}\theta,\qquad\psi=\frac\theta{1-u^2},\qquad c=\tau\psi .
\]
Write $\epsilon=d^4/(\kappa\theta^2)\le1/C_1$. Since $\eta=8\epsilon\theta/d^3$ and $\theta\le d(1-\rho)$, sufficiently small absolute $\epsilon$ gives $\rho<r<1$, $0<\nu<1$, and $\mu=O(\epsilon)$. The bounds $d\eta\le1-r^{2d}\le2d\eta$ give
\begin{equation}\label{eq:ef-bounds}
\psi\ge\frac{d^2}{16\epsilon},\qquad\eta\psi\le\frac1d,\qquad\tau\psi\eta\ge2,\qquad\mu\tau=\frac\eta2 .
\end{equation}
All $O(\epsilon)$ bounds below have absolute constants, uniformly in $d\ge3$ and $0<\theta<1$.

\emph{The cyclic tail.} Let $e_0,\dots,e_{d-1},b_0,\dots,b_{d-1}$ be an orthonormal basis of $\R^{2d}$, let $Q$ be the orthogonal map with $Qe_j=e_{j+1}$, $Qb_j=b_{j+1}$ for $j<d-1$, $Qe_{d-1}=e_0$ and $Qb_{d-1}=-b_0$, and let $v=\sqrt\alpha\,e_0+\sqrt{1-\alpha}\,b_0$. Then
\begin{equation}\label{eq:ef-corr}
\ip{v}{Q^\ell v}=0\quad\text{if }d\nmid\ell,\qquad\ip{v}{Q^{sd}v}=\alpha+(1-\alpha)(-1)^s .
\end{equation}
The tail gradients $\hat g_j=r^jQ^jv$, $j\in\mathbb Z$, are therefore orthogonal at the lags $1,\dots,d-1$. Since $\rho<r$, the stationary filter
\[
Z_k=\sum_{i\le k}\rho^{k-i}\hat g_i
\]
converges. Only the indices $i=k+1-sd$, $s\ge1$, contribute to $\ip{\hat g_{k+1}}{Z_k}$. By the choice of $\alpha$,
\begin{equation}\label{eq:ef-orth}
\begin{aligned}
\ip{\hat g_{k+1}}{Z_k}
&=\frac{r^{2k+2}}\rho\sum_{s\ge1}\nu^s\bigl(\alpha+(1-\alpha)(-1)^s\bigr)\\
&=\frac{r^{2k+2}}\rho\Bigl(\frac{\alpha\nu}{1-\nu}-\frac{(1-\alpha)\nu}{1+\nu}\Bigr)=0 .
\end{aligned}
\end{equation}
Let $w=(I-rQ)^{-1}v=\sum_{j\ge0}r^jQ^jv$, $\hat x_k=\tau r^kQ^kw=\tau\sum_{j\ge k}\hat g_j$ and $\hat f_k=\frac c2r^{2k}$. Then $\hat x_{k+1}-\hat x_k=-\tau\hat g_k$. By \eqref{eq:ef-corr} and $\nu u=\rho^d=1-\theta$,
\[
\ip vw=\frac\alpha{1-u}+\frac{1-\alpha}{1+u}=\frac{1-\nu u}{1-u^2}=\psi,\qquad\norm w^2=\frac{2\psi-1}{1-r^2},
\]
where the second identity follows by squaring $w=v+rQw$. Hence $\ip{\hat g_k}{\hat x_k}=cr^{2k}=2\hat f_k$.

\emph{Residuals of the tail.} The tail is an orbit \eqref{eq:orbit} with $X=\tau w$ and $G=v$, and $Q^{2d}=I$, so $K_\ell^\pm$ depends only on $j=\ell\bmod2d$, and we write $K_j^\pm$. Since $(1-r^\ell)^2>0$, by \eqref{eq:orbit-res} all residuals among the tail triples and the minimizer are positive once $C_\star>0$ and $K_j^\pm>0$ for $1\le j<2d$. Here
\[
\frac{C_\star}{c}=\frac12-\frac{c^{-1}+(2-\psi^{-1})/[2(2-\eta)]-2\mu}{2(1-\mu)}=\frac14+O(\epsilon)>0.
\]
The quadratic part of $K_j^\pm$ is
\[
P_j=\frac{\norm{v-Q^jv}^2+\mu\tau^2\norm{w-Q^jw}^2-2\mu\tau\ip{v-Q^jv}{w-Q^jw}}{2(1-\mu)}.
\]
Since $(I-Q)w=v-\eta Qw$ and $\eta^2\norm w^2\le\eta\psi\le1/d$, we have $\norm{w-Q^jw}\le4d$ for $1\le j<2d$. Thus $\mu\tau^2\norm{w-Q^jw}^2=O(\epsilon)$ and the mixed term is $O(d\eta)=O(\epsilon)$. As $\norm{v-Q^jv}^2=2$ for $j\ne d$ and $\norm{v-Q^dv}^2\le4$, we obtain $P_j\le1+O(\epsilon)$ for $j\ne d$ and $P_d\le2+O(\epsilon)$. Moreover
\[
K_j^+=\tau\bigl(\psi-\ip{Q^jv}w\bigr)-P_j,\qquad K_j^-=\tau\bigl(\psi-\ip v{Q^jw}\bigr)-P_j .
\]
Since $w=v+rQw$, the numbers $a_j=\ip{Q^jv}w$ satisfy $a_0=\psi$ and $a_j=\ip{Q^jv}v+ra_{j-1}$. By \eqref{eq:ef-corr} and $\ip v{Q^dv}=-\nu$ this gives $a_j=r^j\psi$ for $j<d$ and $a_j=r^j\psi-r^{j-d}\nu$ for $d\le j<2d$, so $a_j\le r^j\psi$, and $\ip v{Q^jw}=a_{2d-j}$. Both inner products are therefore at most $r\psi$ for $j\ne d$ and at most $u\psi$ for $j=d$, and \eqref{eq:ef-bounds} and $\psi(1-u)=\theta/(1+u)$ give
\[
K_j^\pm\ge1-O(\epsilon)>0\quad(j\ne d),\qquad
K_d^\pm\ge\frac{\tau\theta}{1+u}-P_d\ge2d-2-O(\epsilon)>0.
\]

\emph{Cold start.} Let $E_0,E_1,E_2\perp\R^{2d}$ be orthonormal, and $G_i=10^{-i/2}E_i$. For $\vartheta>0$ put
\[
g_0=\vartheta^{-1}G_0,\qquad g_1=Z_1+\vartheta^{-1}G_1,\qquad g_2=\hat g_2+\vartheta^{-1}G_2,\qquad g_k=\hat g_k\quad(k\ge3),
\]
$x_k=\hat x_k$ and $f_k=\hat f_k$ for $k\ge3$, and $x_i=x_{i+1}+\tau g_i$ and $f_i=\hat f_3+\frac\tau{2\vartheta^2}\sum_{h=i}^210^{-h}$ for $i=2,1,0$. Then $x_{k+1}-x_k=-\tau g_k$ for every $k\ge0$, and the data are bounded. The map $(x,g,f)\mapsto(\vartheta x,\vartheta g,\vartheta^2f)$ multiplies all residuals by $\vartheta^2$. As $\vartheta\to0$ it sends the tail triples to $(0,0,0)$ uniformly in $k$ and the first three triples to $X_i=\tau\sum_{h=i}^2G_h$, $G_i$, $F_i=\frac\tau2\sum_{h=i}^210^{-h}$. Put $(X_3,G_3,F_3)=(0,0,0)$. These limiting data are obtained from Lemma~\ref{lem:prefix} with $t_i=2$ and strong convexity parameter $\mu'=\mu\tau/2$, by multiplying the points and values by $\tau/2$. The assumptions hold because $\tau\ge12$ and $\mu\tau=O(\epsilon)\le1/6$ for small $\epsilon$. The residuals for $\Fc_{\mu,2/\tau}$ are $\tau/2$ times those of the lemma and do not decrease when the smoothness bound is increased to $1$. By continuity and uniform convergence, for small $\vartheta$ all off-diagonal residuals that involve one of the first three triples are positive. The residuals among the tail triples and the minimizer do not depend on $\vartheta$. Lemma~\ref{lem:interp} gives $F\in\Fc_{1/\kappa,1}(\R^{2d+3})$ with minimizer $0$ that interpolates the data.

\emph{The trajectory.} Every step is a multiple of the current gradient, so the window at iteration $k$ is spanned by $g_{k-m},\dots,g_k$, whose lags to $k+1$ lie between $1$ and $d-1$. By Lemma~\ref{lem:transfer} it suffices to show that $g_{k+1}$ is orthogonal to these gradients and to $z_k$ and $y_k$. Among the tail gradients this is \eqref{eq:ef-corr}. The cold gradients are mutually orthogonal, since $\ip{Z_1}{\hat g_2}=0$ by \eqref{eq:ef-orth}. For $3\le j\le d$ we have $\hat g_j\perp Z_1$, because $Z_{j-1}=\rho^{j-2}Z_1+\sum_{i=2}^{j-1}\rho^{j-1-i}\hat g_i$, $\hat g_j\perp Z_{j-1}$ and $\hat g_j\perp\hat g_i$ for these lags. The registers satisfy $z_1=g_1$, $y_1=-\rho\tau g_0$, $z_2=\rho z_1+g_2$ and $y_2=-\tau\rho(\rho g_0+g_1)$, so by induction their components in $\R^{2d}$ are $Z_k$ and $-\tau\rho Z_{k-1}$ for $k\ge2$. For $k\ge2$, $\hat g_{k+1}$ is orthogonal to $Z_k$ by \eqref{eq:ef-orth} and to $Z_{k-1}=(Z_k-\hat g_k)/\rho$. The components of the registers along $E_0,E_1,E_2$ are orthogonal to the tail gradients, and for $k\le1$ the registers lie in the span of $g_0$ and $g_1$. Hence $g_{k+1}\perp S_k$ for every admissible rule, and every such rule generates the trajectory. With the constructed $F$ fixed, $F^\star=0$ and $F(x_k)-F^\star=\frac c2r^{2k}$ for $k\ge3$, so the fixed initial segment again leaves the asymptotic coefficient unchanged. Since $0<\eta\le1/2$ for sufficiently large $C_1$ and $-\log(1-\eta)\le2\eta$, we obtain $\lim_{\varepsilon\downarrow0}N(\varepsilon)/\log(1/\varepsilon)=1/(-2\log r)\ge1/(4\eta)=\kappa\theta/(32d)$.
\end{proof}

\subsection{The right time scale}\label{sec:lower}

Theorem~\ref{thm:slow} excludes memories that forget too slowly and Theorem~\ref{thm:fast} a pair that forgets too fast, which leaves a single time scale.

\begin{corollary}[long memory or fast forgetting]\label{cor:memA}
For every $\kappa$ let a method use one memory of Definition~\ref{def:separable} with window $m=m(\kappa)$ and forgetting index $\delta=\delta(\kappa)$, and let $\limsup_{\varepsilon\downarrow0}N(\varepsilon)/\log(1/\varepsilon)\le C\sqrt\kappa$ on all of $\Fc_{1/\kappa,1}$ in every dimension, with $C$ independent of $\kappa$. Then, with constants depending only on $C$ and for all large $\kappa$:
\begin{enumerate}
\item[(a)] $m\ge c\sqrt\kappa$ or $\delta\ge c/\sqrt\kappa$.
\item[(b)] If the memory consists of the window and one pair as in Theorem~\ref{thm:fast} and $m=o(\kappa^{1/4})$, then $\Omega(\kappa^{-1/2})\le1-\rho\le O((m+1)\kappa^{-1/2})$. In particular $1-\rho=\Theta(\kappa^{-1/2})$ for a bounded window.
\end{enumerate}
\end{corollary}
\begin{proof}
(a) Choose $p=\max\{8,m+2,\lfloor\kappa\delta\rfloor+1\}$, so $p\le\max\{m+8,\kappa\delta+1\}$ and $\delta<p/\kappa$. If $\kappa<C_0p^2$, then $p\gtrsim\sqrt\kappa$ directly; otherwise Theorem~\ref{thm:slow} and the assumed accelerated bound give $\kappa/p\lesssim C\sqrt\kappa$. Thus $p\gtrsim\sqrt\kappa$ in either case, which proves (a) for all large $\kappa$.

(b) The pair has normalized weights $\rho^{-i}$, so $\delta=1-\rho$, and (a) gives $\delta=\Omega(\kappa^{-1/2})$ since $m=o(\sqrt\kappa)$. Enlarge a window of length zero to one if necessary. Put $d=m+2$ and $\theta=1-\rho^d$. Either $\kappa\theta^2<C_1d^4$, or Theorem~\ref{thm:fast} gives $\kappa\theta/d\lesssim C\sqrt\kappa$. Hence $\theta=O(d^2/\sqrt\kappa)=o(1)$ under the stated assumption. Therefore
\[
\delta=1-(1-\theta)^{1/d}\le\frac{-\log(1-\theta)}d=O\left(\frac d{\sqrt\kappa}\right).
\tag*{\qed}
\]
\let\qed\relax
\end{proof}

By \eqref{eq:accel} every uniformly accelerated method satisfies the assumption with $C=1/c$, so the corollary applies to it. Both alternatives can be attained by accelerated constructions. A window of length $1.39\sqrt\kappa$ is accelerated (Proposition~\ref{prop:window}), and so is a filter pair with $1-\rho\in[1/(2\sqrt\kappa),1/\sqrt\kappa]$ (Lemma~\ref{lem:pair}, whose proof uses only $1-\rho\le1/\sqrt\kappa$). A register that the method never uses increases $\delta$ without changing the iterates, so part (a) is only a necessary condition. The activation schedule and the smooth pair of Section~\ref{sec:filters} lie outside Definition~\ref{def:separable}: the first has unboundedly many registers, and the second adapts its weights to the run. For memories whose factors change in time, $\delta$ gives no upper bound like part (b).

\begin{remark}[periodic attenuation]\label{rem:attenuation}
There are absolute constants $A>0$, $\eta\in(0,1)$ and $c,C>0$ with the following property. Let $\kappa\ge4$, $h=\lceil\sqrt\kappa\,\rceil$ and $\rho=1-1/h$, let $S_k$ contain $\nabla f(x_k)$ and the pair of Theorem~\ref{thm:fast}, and multiply both registers by $\beta=\eta/h$ every $T=\lceil Ah\rceil+1$ iterations. This memory is fixed in advance but depends on $\kappa$, and every such method satisfies $f(x_k)-f^\star\le C e^{-ck/\sqrt\kappa}(f(x_0)-f^\star)$ on $\Fc_{\mu,L}$ with $L/\mu=\kappa$ (Section~\ref{sec:filters}), although $\delta=1-\rho\beta^{1/T}\asymp\log\kappa/\sqrt\kappa$.
\end{remark}

In our lower bounds the memory, including the number of aggregates, is fixed before the instance. Since $n$ gradient sums with weights $i^a$, $a<n$, terminate on $\R^n$, families whose number of aggregates grows with the dimension are not covered. One function serves all accuracies, so no prefactor depending on $\kappa$ compensates the bound. The convex rate of SESOP (Proposition~\ref{prop:sesop}) does not exclude such a slow tail, and the contraction $(1-3/\kappa)^2$ on the instances of Corollary~\ref{cor:sesop} is of the order of the worst case $1-4/\kappa+O(\kappa^{-2})$ of GD\nobreakdash-ELS \cite{deklerk2017worst}. Like that worst case, our bounds concern the asymptotic rate. The constructions first reduce the gap by a factor that depends on $\kappa$, so they say nothing about moderate accuracy. Section~\ref{sec:smooth} extends the obstruction to infinitely differentiable objectives.

\subsection{Smoothness does not help}\label{sec:smooth}

The interpolation constructions in Theorems~\ref{thm:slow} and~\ref{thm:fast} guarantee $C^{1,1}$ regularity but do not establish $C^2$ regularity near the minimizer. Near a minimizer with a continuous Hessian a function looks like a quadratic, on which SESOP coincides with conjugate gradients, but the memory accumulated before the quadratic region keeps SESOP slow. Both instances below rest on one exact plane step on a quadratic with two eigenvalues.
\begin{lemma}[plane rule on two eigenvalues]\label{lem:invariant}
Let $B=I+bvv^\top$ on $\R^3$ with $\norm v=1$ and $b>0$, and let $q(x)=\frac12x^\top Bx$. Let $g=Bx$ and the inherited step $s$ be nonzero and orthogonal. Put $\gamma=\ip vg^2/\norm g^2$ and $\zeta=\ip vs^2/\norm s^2$, and assume $\gamma,\zeta>0$ and $\gamma+\zeta<1$. Let $x^+$ minimize $q$ over $x+\spn\{g,s\}$, with gradient $g^+$ and step $s^+=x^+-x$. Then
\[
\frac{\ip v{g^+}^2}{\norm{g^+}^2}=1-\gamma-\zeta,\qquad
\frac{\ip v{s^+}^2}{\norm{s^+}^2}=\frac\gamma{(1+b\zeta)^2+b^2\gamma\zeta}.
\]
The quantity $J=b\gamma\zeta(1-\gamma-\zeta)/(1+b\zeta)^2$ is invariant under this update, and
\[
\frac{q(x^+)}{q(x)}=\frac{b^2\gamma(1-\gamma-\zeta)}{(1+b-b\gamma)(1+b\gamma+b\zeta)}\ge\frac{bJ}{1+bJ}.
\]
If $1+b=K\ge12$ and $r\in(0,1)$ solves $K(1-r)=(1+r)^2$, then $\gamma=\frac{1+r}{3+r}$, $\zeta=\frac{1-r}{3+r}$ is a fixed point, at which $q(x^+)/q(x)=r^2$ and $J=r/(3+r)^2>\frac1{32}$.
\end{lemma}

\begin{proof}
The gradient $g^+$ is orthogonal to $g$ and $s$, so it is parallel to the unit normal $h$ of the plane, with $\ip vh^2=1-\gamma-\zeta$. Write $s^+=-\alpha g+\beta s$. Then $g^+=g+Bs^+\perp s$ gives
\[
\beta\norm s^2(1+b\zeta)=\alpha b\ip vg\ip vs,\qquad\text{hence}\qquad\ip v{s^+}=-\frac{\alpha\ip vg}{1+b\zeta},
\]
and $\norm{s^+}^2=\alpha^2\norm g^2+\beta^2\norm s^2$ gives the second update formula. Substitution of the updated ratios into $J$ proves its invariance. Since $B^{-1}=I-\frac b{1+b}vv^\top$ and $B^{-1}(g^+-g)=s^+$ is orthogonal to $h$, we get $g^+=-\frac{b\ip vg\ip vh}{1+b-b\ip vh^2}h$, and $q=\frac12g^\top B^{-1}g$ gives the ratio. The identity $(1+b-b\gamma)(1+b\gamma+b\zeta)=b^2\gamma(1-\gamma-\zeta)+(1+b)(1+b\zeta)$ and $(1+b)\zeta\le1+b\zeta$ give the lower bound. For the fixed point, $b(1-r)=r(3+r)$ gives $b\zeta=r$ and $b\gamma=r(1+r)/(1-r)$, and the formulas follow. The function $r/(3+r)^2$ increases on $(0,1)$, and $K\ge12$ gives $r>\frac12$.
\end{proof}

With $b=\kappa-1$ the tail of the following explicit instance for the plane rule sits at this fixed point, and its local and global condition numbers are comparable.

\begin{proposition}[plane rule in $\R^3$]\label{prop:plane3}
Let $\kappa\ge100$, $q(x)=\frac12x^\top Dx$ with $D=\operatorname{diag}(1,1,\kappa)$, and let $r\in(0,1)$ solve $(1+r)^2=\kappa(1-r)$. There are a convex function $\phi\in C^\infty(\R^3)$, $\phi\ge0$, that vanishes on a neighbourhood of the unit ball, and a point $x_0$ such that $f=q+\phi$ is $1$-strongly convex with a $(\kappa+C\sqrt\kappa)$-Lipschitz gradient and a $C\sqrt\kappa$-Lipschitz Hessian for an absolute constant $C$, and the plane rule on $f$ started at $x_0$ satisfies $f(x_{k+1})=r^2f(x_k)$ for $k\ge1$. Here $r^2=1-8/\kappa+O(\kappa^{-2})$ and $\liminf_{\varepsilon\downarrow0}N(\varepsilon)/\log(1/\varepsilon)\ge(\kappa-4)/8$.
\end{proposition}

\begin{proof}
We first construct a slow orbit on the quadratic $q$ with a prescribed inherited step. A convex perturbation outside the unit ball then sends the first line search onto this orbit with the required step direction. Inside the ball the objective is exactly quadratic, but the inherited direction preserves the slow orbit.

\emph{The quadratic orbit.} Let $\cos\alpha=(1+r)/2$ and $\eta=\sqrt{\cos\alpha}/\kappa$. Set
\[
\hat x_k=(r^k\cos k\alpha,\ r^k\sin k\alpha,\ (-r)^k\eta),\qquad k\in\mathbb Z.
\]

Put $\hat g_k=D\hat x_k$ and $\hat s_k=\hat x_k-\hat x_{k-1}$. Direct substitution, using $\kappa(1-r)=(1+r)^2$, gives
\[
\hat s_{k+1}=-(1-r)\hat g_k+r^2\hat s_k,\qquad
\ip{\hat g_{k+1}}{\hat g_k}=\ip{\hat g_{k+1}}{\hat s_{k+1}}=0.
\]
Hence $\hat g_{k+1}\perp\hat s_k$, so $\hat x_{k+1}$ minimizes $q$ over $\hat x_k+\spn\{\hat g_k,\hat s_k\}$, and the orbit formula gives $q(\hat x_{k+1})=r^2q(\hat x_k)$.

\emph{Entering the orbit from a cold start.} We need a nonnegative convex perturbation that vanishes near the unit ball and directs the first line search along $\hat s_1$. Put $d=\hat s_1/\norm{\hat s_1}$ and $x_0=\hat x_1-8d$. Directly from the orbit,
\[
\norm{\hat x_1}\le1,\qquad\norm{D\hat x_1}\le\sqrt2,\qquad\norm{Dd}\le\sqrt{\kappa+1}.
\]
Choose $\lambda=A\sqrt\kappa$ with a sufficiently large absolute $A$ and set $w_0=-8\lambda d-Dx_0$. Since $\norm{Dx_0}=O(\sqrt\kappa)$, the unit vector $u=w_0/\norm{w_0}$ is uniformly close to $-d$ as $A$ increases. Thus $\ip u{x_0}>3$ and $\norm{w_0}=O(\sqrt\kappa)$, with absolute constants.

Fix a smooth convex function $\psi$ with $\psi=0$ on $(-\infty,0]$ and $\psi'=1$ on $[1,\infty)$; one obtains it by integrating twice a nonnegative smooth function supported in $(0,1)$ with integral one. Then
\[
\phi(x)=\norm{w_0}\psi(\ip ux-2)
\]
is nonnegative, convex and $C^\infty$, vanishes on $\{\norm x<2\}$, and satisfies $\nabla\phi(x_0)=w_0$. Its Hessian norm and Hessian Lipschitz constant are $O(\sqrt\kappa)$ because $\psi''$ and $\psi^{(3)}$ are bounded. Hence $\nabla(q+\phi)(x_0)=-8\lambda d$, with the required global regularity bounds.

\emph{Verification of the trajectory and rate.} Since $f\ge q$ and $f=q$ at every $\hat x_k$ with $k\ge1$, each prescribed quadratic minimizer also minimizes $f$ on the same affine space. The first line search reaches $\hat x_1=x_0+8d$, because $\ip{D\hat x_1}d=0$, and produces $s_1\parallel\hat s_1$. Induction therefore gives $x_k=\hat x_k$ for $k\ge1$. Finally $f^\star=0$ and $-\log r\le(1-r)/r\le4/(\kappa-4)$, so $\lim_{\varepsilon\downarrow0}N(\varepsilon)/\log(1/\varepsilon)=1/(-2\log r)\ge(\kappa-4)/8$.
\end{proof}

For SESOP, three steps bring the state near the fixed point of Lemma~\ref{lem:invariant}.

\begin{proposition}[smooth instances for SESOP]\label{prop:smooth}
Let SESOP use the previous step, the anchor and one aggregate with slowly forgetting weights such that $w(0)\ge0$ and $w(1),w(2)>0$. There is $\mu_w>0$ depending only on $w(0),w(1),w(2)$ such that for every $K\ge12$ there are $G\in C^\infty(\R^5)$ and a starting point with the following properties. The Hessian satisfies $\mu_wI\preceq\nabla^2G\preceq2KI$, all derivatives of order at least two are bounded, and $G(x)=\frac12x^\top Hx$ near the minimizer $0$ with $H=I+(K-1)ee^\top$, $\norm e=1$. From $x_3$ on SESOP stays in this region, and $G(x_{k+1})\ge\frac{K-1}{K+127}G(x_k)$ for $k\ge3$.
\end{proposition}
\begin{proof}
The anchor and aggregate remove two directions from the five-dimensional search geometry, leaving a three-dimensional plane method. Its subspace changes with the memory, so the key point is stability: small total motion of the subspace keeps the invariant of Lemma~\ref{lem:invariant} positive and hence keeps convergence slow. Steps 1 and 2 give a sufficient condition on the state at $x_3$. Steps 3 and 4 construct three steps satisfying it, and Step 5 realizes those steps on a smooth function that is exactly quadratic near the tail.

\emph{Step 1: eliminate the stored directions.} From now on SESOP has one aggregate $\sigma_k=\sum_{i\le k}w(i)g_i$, so $S_k=\spn\{g_k,s_k,x_k-x_0,\sigma_{k-1}\}$, since $\sigma_k=\sigma_{k-1}+w(k)g_k$, and $g_k\perp s_k,\,x_k-x_0,\,\sigma_{k-1}$. Let $H=I+aee^\top$ on $\R^5$ with $a=K-1$, and $f(x)=\frac12x^\top Hx=\frac12\norm y^2$ in energy coordinates $y=H^{1/2}x$. Put
\begin{gather*}
C_k=\bigl[y_k-H^{1/2}x_0,\ H^{1/2}\sigma_{k-1}\bigr],\qquad W_k=(\operatorname{range}C_k)^\perp,\\
P_k=P_{W_k},\qquad\lambda_k=\norm{P_ke}^2,\qquad t_k=P_k(y_k-y_{k-1}).
\end{gather*}
While the objective equals $f$, the orthogonality relations give $y_k\in W_k$ and $y_k\perp t_k$. The search directions in energy coordinates are $Hy_k$, $y_k-y_{k-1}$ and the columns of $C_k$, so $y_{k+1}$ minimizes $\frac12\norm y^2$ over $y_k+\spn\{B_ky_k,t_k\}$ with $B_k=P_kH|_{W_k}=I+a(P_ke)(P_ke)^\top$. If $C_k$ has rank two, then $\dim W_k=3$, and in the coordinates $B_k^{-1/2}y$ this is the plane step of Lemma~\ref{lem:invariant} with $b=a\lambda_k$, $v=P_ke/\norm{P_ke}$, gradient $B_k^{1/2}y_k$ and inherited step $B_k^{-1/2}t_k$. Define $J(P_k,y_k,t_k)$ by evaluating the invariant in Lemma~\ref{lem:invariant} at this transformed gradient and step, and write $J_k=J(P_k,y_k,t_k)$. The step conserves $J(P_k,\cdot,\cdot)$, so $J_{k+1}$ differs from it only through the projector. To extend $J$ off projectors, use $B_M=I+a(Me)(Me)^\top$, $v=Me/\norm{Me}$, $b=a\norm{Me}^2$ and the vectors $B_M^{1/2}y$, $B_M^{-1/2}Mt$ in the same algebraic formula. For nonzero $y,t$, this extension is smooth and homogeneous of degree zero in $y,t$ where $\norm{Me}\ge\frac14$ and $\norm{Mt}\ge\frac34\norm t$. These conditions hold for $t\in\operatorname{range}P$ on the segment between projectors $P,P'$ with $\norm{P'-P}\le\frac14$ and $\norm{Pe}^2\ge\frac14$, so by compactness there is $L_K$, depending only on $K$, such that
\begin{equation}\label{eq:J-lip}
\bigl|J(P',y,P't)-J(P,y,t)\bigr|\le L_K\norm{P'-P}
\end{equation}
for such $P,P'$ and nonzero vectors $y\in\operatorname{range}P\cap\operatorname{range}P'$ and $t\in\operatorname{range}P$ with $t\perp y$.

\emph{Step 2: a sufficient condition for a persistent slow tail.} The stored columns must remain independent, retain the high-curvature direction in their orthogonal complement, and move little enough to keep $J_k$ positive. These requirements become quantitative in \eqref{eq:basin}. Let $q_K=\sqrt{1-1/K}$ and $E=\sum_{j\ge0}(2+K|w(3+j)|)q_K^j$, which is finite for slowly forgetting weights. Let the objective be convex and equal to $f$ near the ball $\{\norm{H^{1/2}x}\le R\}$, where $R=\norm{y_3}$, let $g_3\perp s_3,\,x_3-x_0,\,\sigma_2$ and $t_3\ne0$, and assume that for some $\sigma_*>0$
\begin{equation}\label{eq:basin}
\sigma_{\min}(C_3)>\sigma_*,\qquad\lambda_3>\tfrac34,\qquad J_3>\tfrac1{32},\qquad RE<\min\{\sigma_*/8,\,\sigma_*/(128L_K)\}.
\end{equation}
We show that then SESOP never terminates and $f(x_{k+1})\ge\frac{K-1}{K+127}f(x_k)$ for $k\ge3$. The minimizer of $f$ over $x_k+S_k$ is not worse than the gradient step $x_k-g_k/K$, so its energy norm is at most $q_K\norm{y_k}$. Hence it lies where the objective equals $f$, the gradient of the objective there is orthogonal to $S_k$, and by convexity it is the SESOP iterate. By induction $\norm{y_{3+j}}\le q_K^jR$. The columns of $C_k$ change by $y_{k+1}-y_k$ and $w(k)Hy_k$, so $\sum_{k\ge3}\norm{C_{k+1}-C_k}\le RE<\sigma_*/8$, and $\sigma_{\min}(C_k)>\sigma_*/2$ for all $k$. Since $C_k$ and $C_{k+1}$ both have rank two, the projector identity gives
\[
\norm{P_{k+1}-P_k}=\norm{P_{k+1}(C_k-C_{k+1})C_k^+}
\le\frac{2}{\sigma_*}\norm{C_{k+1}-C_k},
\]
where $C_k^+=(C_k^\top C_k)^{-1}C_k^\top$ is the left pseudoinverse. Hence $\sum_{k\ge3}\norm{P_{k+1}-P_k}<\frac14$ and $\lambda_k>\lambda_3-\frac14>\frac12$. To apply \eqref{eq:J-lip}, take $P=P_k$, $P'=P_{k+1}$, $y=y_{k+1}$ and $t=y_{k+1}-y_k$. Here $y_k,y_{k+1}\in W_k$, $y_{k+1}\in W_{k+1}$ and $y_{k+1}\perp t$, while $P_{k+1}t=t_{k+1}$. By Lemma~\ref{lem:invariant} and \eqref{eq:J-lip} the invariant changes by less than $L_K\cdot\frac{2}{\sigma_*}RE<\frac1{64}$ in total, so $J_k>\frac1{64}$. A positive invariant excludes a vanishing next point, and a projector change below $\frac14$ excludes a vanishing projected step, which closes the induction. Finally $b=a\lambda_k\ge a/2$ and Lemma~\ref{lem:invariant} give
\[
\frac{f(x_{k+1})}{f(x_k)}\ge\frac{bJ_k}{1+bJ_k}\ge\frac{a/128}{1+a/128}=\frac{K-1}{K+127}.
\]

\emph{Step 3: choose the limiting geometry.} We place the anchor and aggregate in two independent directions orthogonal to the high-curvature direction, and place the remaining state at the fixed point of Lemma~\ref{lem:invariant}. Write $w_i=w(i)$ and $c=(w_0+w_2)/w_1>0$. Choose $\epsilon_g>0$ sufficiently small depending only on $w_0,w_1,w_2$, as specified in Step~5, and set $\rho_w=\epsilon_g^{3/2}$ and $\mu_w=\epsilon_g^2<1/2$. Let $e_0,e_1,e_2,v_1,v_2$ be an orthonormal basis of $\R^5$, $d^\ast=(-e_0+e_1-e_2)/\sqrt3$, $d=2$ if $w_0\ne2w_2$ and $d=3$ otherwise, and $U=(de_0+(d+1)e_1+e_2)/\sqrt3$. With the fixed point $(\gamma^\ast,\zeta^\ast)$ of Lemma~\ref{lem:invariant} for $b=a$ let $e=\sqrt{\gamma^\ast}\,v_1+\sqrt{\zeta^\ast}\,d^\ast+\sqrt{\gamma^\ast}\,v_2$ and $H=I+aee^\top$. The vectors $U$ and $V=w_0e_0+(w_0+w_2)e_1+w_2e_2$ are orthogonal to $d^\ast$ and $e$, and independent since $d\ne w_0/w_2$.

\emph{Step 4: reach the required state in three steps.} The parameter $\theta$ makes $x_3$ small while the stored columns stay independent. The projections below enforce exact subspace optimality without changing the limiting geometry. For small $\theta>0$ put $n=U/\norm U$ and choose the small root $\nu=O(\theta^2)$ of $\nu(\norm U-\nu)=\theta^2\ip{v_1}{H^{-1}v_1}$. Set
\[
x_3=\theta H^{-1}v_1+\nu n,\qquad g_3=Hx_3=\theta v_1+\nu n.
\]
Then $\ip{g_3}{d^\ast}=\ip{g_3}{x_3-U}=0$. Put $p=x_3-d^\ast-U$, let $t_\theta$ be the projection of $e_2$ onto $\spn\{g_3,p,d^\ast\}^\perp$, $E_\theta=\spn\{p,d^\ast,t_\theta\}$, and
\[
g_2=P_{E_\theta\cap p^\perp}(\epsilon_ge_2),\qquad g_0=\epsilon_g\frac{P_{E_\theta\cap g_2^\perp}e_0}{\norm{P_{E_\theta\cap g_2^\perp}e_0}},\qquad g_1=\epsilon_gc\,\frac{P_{E_\theta\cap\{g_0,g_2\}^\perp}e_1}{\norm{P_{E_\theta\cap\{g_0,g_2\}^\perp}e_1}},
\]
and $x_0=U$, $x_1=U-\frac{d-1}{\sqrt3\,\epsilon_g}g_0$, $x_2=x_3-d^\ast$. As $\theta\downarrow0$, $x_3\to0$, $E_\theta\to\spn\{e_0,e_1,e_2\}$, $(g_0,g_1,g_2)\to\epsilon_g(e_0,ce_1,e_2)$, $x_1\to X_1=(e_0+(d+1)e_1+e_2)/\sqrt3$ and $x_2\to X_2=-d^\ast$, and for small $\theta$ the vectors $g_0,g_1$ are independent and so are $p,x_2-x_1,g_2$. Let a strongly convex function have the gradients $g_i$ at $x_i$ for $i\le2$ and equal $f$ near $x_3$. In each of the first three SESOP steps the new point lies in the search space and its gradient is orthogonal to it:
\begin{align*}
&S_0=\spn\{g_0\}, && g_1\perp g_0,\\
&S_1=\spn\{g_0,g_1\}\ni x_2-x_1=p+\tfrac{d-1}{\sqrt3\,\epsilon_g}g_0, && g_2\perp g_0,g_1,\\
&S_2=\spn\{g_2,x_2-x_1,p,\sigma_1\}=E_\theta\ni d^\ast, && g_3\perp d^\ast,p,t_\theta .
\end{align*}
Hence SESOP generates $x_1,x_2,x_3$. As $\theta\downarrow0$, $C_3\to[-U,\epsilon_gV]$, whose range is orthogonal to $e$, so $\lambda_3\to1$. Since $\nu=O(\theta^2)$ and $H^{1/2}$ fixes $e^\perp$, the compressed gradient and inherited step at $k=3$ become parallel to $v_1$ and $d^\ast$, whose coordinates are $(\gamma^\ast,\zeta^\ast)$, so $J_3\to r/(3+r)^2>\frac1{32}$, and $R=\norm{y_3}\to0$. With $\sigma_*=\frac14\sigma_{\min}[-U,\epsilon_gV]$ all conditions of Step~2 hold for small $\theta$, provided that the objective equals $f$ on a ball of fixed radius around $0$.

\emph{Step 5: realize the finite data and preserve the quadratic tail.} We interpolate with strictly separated conjugate branches, making the function locally quadratic at the prescribed points and near $0$, so that smoothing preserves their gradients and the tail. Let $L=2K$, $\Lambda=L-\mu_w$, $(F_0,F_1,F_2)=(3,2,1)\rho_w$, $X_0=U$ and $G_i=\epsilon_g(e_0,ce_1,e_2)_i$. With $q_i=g_i-\mu_wx_i$ and $u_i=\ip{x_i}{q_i}-F_i+\frac{\mu_w}2\norm{x_i}^2$ put
\[
\psi(\xi)=\max\Bigl\{\tfrac12\xi^\top(H-\mu_wI)^{-1}\xi,\ \max_{i\le2}\Bigl[u_i+\ip{x_i}{\xi-q_i}+\frac{\norm{\xi-q_i}^2}{2\Lambda}\Bigr]\Bigr\}
\]
and $F(x)=\psi^\ast(x)+\frac{\mu_w}2\norm x^2$. The limit points $X_i$ are bounded and $\norm{G_i}=O_w(\epsilon_g)$, where $O_w$ depends only on $w_0,w_1,w_2$. Since $L\ge24$ and $\norm{(H-\mu_wI)^{-1}}\le(1-\mu_w)^{-1}$, all quadratic correction terms below are $O_w(\epsilon_g^2)$ uniformly in $K\ge12$. For $i<j$, $\ip{G_j}{X_i-X_j}=0$, so $C_{ij}=(j-i)\rho_w+O_w(\epsilon_g^2)>0$. For $i>j$,
\[
-\ip{G_j}{X_i-X_j}\ge\frac{\epsilon_g\min\{1,c\}}{\sqrt3},\qquad
\ip{X_i}{G_i}\ge\frac{\epsilon_g\min\{1,c\}}{\sqrt3}.
\]
The first inequality gives $C_{ij}>0$ for $i>j$. The second holds for every $0\le i\le2$ and makes the $i$th data branch dominate the quadratic branch at $q_i$: the leading positive term is of order $\epsilon_g$, whereas $F_i=O(\rho_w)$ and the remaining terms are $O_w(\epsilon_g^2)$. At $0$ each data branch equals $-F_i+O_w(\epsilon_g^2)<0$. These signs follow from $\epsilon_g^2=o(\rho_w)$ and $\rho_w=o(\epsilon_g)$, so one choice of $\epsilon_g$ works for every $K\ge12$. Fix that choice. The difference between data branches $i$ and $j$ at $q_i$ is $C_{ji}$, so these inequalities give strict separation, which persists for the actual data when $\theta$ is sufficiently small.

Every branch of $\psi$ is $\Lambda^{-1}$-strongly convex, so $F$ is $\mu_w$-strongly convex with an $L$-Lipschitz gradient. Strict separation gives
\[
F(x_i+h)=F_i+\ip{g_i}h+\frac L2\norm h^2\quad(i\le2),\qquad F(x)=f(x)\quad\text{near }0.
\]
For fixed $K,w$, continuity of the data and the strict margins at $\theta=0$ give positive radii for these neighbourhoods uniform over all sufficiently small $\theta$. Choose an even mollifier $\eta_{\epsilon_s}$ with support smaller than these radii and set
\[
G=F*\eta_{\epsilon_s}-\frac12\int\xi^\top H\xi\,\eta_{\epsilon_s}(\xi)\,d\xi.
\]
Then $\nabla G(x_i)=g_i$, $G=f$ on a fixed neighbourhood of $0$, and $\mu_wI\preceq\nabla^2G\preceq LI$. All higher derivatives are bounded by convolving the bounded weak Hessian of $F$ with derivatives of the kernel. Finally decrease $\theta$ so that the energy ball of radius $R=\norm{y_3}$ lies in this quadratic neighbourhood and the conditions of Step~2 hold. SESOP then generates $x_1,x_2,x_3$ and the slow tail, completing the proof.
\end{proof}

Since $J_3>0$ implies $G(x_3)>0$, Proposition~\ref{prop:smooth} and $\log(1+u)\le u$ give
\[
\liminf_{\varepsilon\downarrow0}\frac{N(\varepsilon)}{\log(1/\varepsilon)}
\ge\frac1{\log((K+127)/(K-1))}\ge\frac{K-1}{128}.
\]
For fixed weights this bound is of the order of the condition number $\kappa\le2K/\mu_w$ of $G$. As $K$ is the condition number near the minimizer, unrestarted SESOP has no local accelerated rate either, unlike the method of \cite{tupitsa2026cacute} with adaptive restarts. The conditions \eqref{eq:basin} are strict inequalities, and an exact step restores the orthogonality relations, so this slow regime survives small perturbations of the state at $x_3$, which the cyclic tails do not.

\section{Memories that accelerate}\label{sec:upper}

We now give memories that hold the virtual momentum of Lemma~\ref{lem:cert}. A window of length of order $\sqrt\kappa$ holds it over blocks of steps (Section~\ref{sec:window}). A pair of geometric filters holds it in two vectors once its factor matches $\kappa$, and Section~\ref{sec:filters} finds this factor by a fixed schedule and by an adapted horizon. Section~\ref{sec:inexact} extends the activation schedule to inexact subproblems.

\subsection{Windows}\label{sec:window}

Corollary~\ref{cor:memA}(a) requires a window of length $\Omega(\sqrt\kappa)$ from a slowly forgetting memory. With $\pi_b=0$ and a constant $a_k=a$, the virtual momentum \eqref{eq:prec} is a combination of the gradients and steps after the origin $b$, so a window that covers a block of steps contains it, and the certificate of Lemma~\ref{lem:cert} can be restarted block by block. Only the length of the window depends on a bound $K\ge\kappa$. We say that $S_k$ \emph{covers the last $m$ gradients and $n$ steps after $b$} if it contains $\nabla f(x_i)$ for $\max(b,k-m+1)\le i\le k$ and $s_i$ for $\max(b+1,k-n+1)\le i\le k$.

\begin{proposition}[long windows suffice]\label{prop:window}
Let the iterates satisfy \eqref{eq:step}, and for some $b\ge0$, $m\ge1$ and every $k\ge b$ let $S_k$ cover the last $m$ gradients and $m$ steps after $b$. If $K\ge\kappa$ and $m\ge\lceil2\log2\,\sqrt K\,\rceil$, then $f(x_{b+h})-f^\star\le2\rho_K^h(f(x_b)-f^\star)$ for $h\ge0$, with $\rho_K=1-1/(2\sqrt K)$.
\end{proposition}
\begin{proof}
For $K=1$ the claim is trivial. Otherwise apply Lemma~\ref{lem:cert} with $\pi_b=0$, $\lambda=L/K\le\mu$ and constant $a_k=a=K^{-1/2}$, so that $a^2/\lambda=1/L$. By \eqref{eq:prec}, $\pi_{b+j}$ is a combination of $\nabla f(x_{b+i})$ and $s_{b+i}$, $1\le i\le j$, and $S_{b+j}$ covers them for $j\le m$. Hence Lemmas~\ref{lem:descent} and~\ref{lem:cert} give, with $B=m+2$,
\[
f(x_{b+h})-f^\star\le\min\bigl\{1,2(1-a)^{h-1}\bigr\}(f(x_b)-f^\star)\qquad(1\le h\le B).
\]
With $\rho=\rho_K$ we have $\rho^2\ge1-a$ and $2\rho^m\le1$, so $f(x_{b+B})-f^\star\le2\rho^{2B-2}(f(x_b)-f^\star)\le\rho^B(f(x_b)-f^\star)$, and $f(x_{b+h})-f^\star\le2\rho^h(f(x_b)-f^\star)$ for $h<B$ by $\rho^2\ge1-a$ if $h\ge2$ and $\rho\ge\frac12$ if $h=1$. Covering after $b'\ge b$ follows from covering after $b$, so chaining blocks of length $B$ gives the claim.
\end{proof}

A window of length $1.39\sqrt\kappa$ therefore suffices (the constant $2\log2\approx1.39$ comes from the crude requirement $2\rho_K^m\le1$), and the order $\Theta(\sqrt\kappa)$ in Corollary~\ref{cor:memA} is sharp for slowly forgetting memories. If the Hessian at the minimizer is a multiple of the identity plus a term of low rank, even a short window is fast near the minimizer, because orthogonal gradients over a block of steps force a quadratic rate over blocks.

\begin{theorem}[local window bound]\label{thm:window-ub}
Let $f\in\Fc_{\mu,L}$, let the iterates satisfy \eqref{eq:step}, and for some $b\ge0$, $m\ge1$ and every $k\ge b$ let $S_k$ cover the last $m$ gradients and $m-1$ steps after $b$. Let $\nabla^2f$ be $M$-Lipschitz on the ball $B=\{x:\norm{x-x^\star}<\omega\}$, and let $\nabla^2f(x^\star)=\chi I+R$ with $\operatorname{rank}R\le m-1$. Then there is a constant $c$, depending only on $m$, $L$, $\mu$, $M$, $\norm R$ and $\omega$, such that $f(x_{k+m})-f^\star\le c\,(f(x_k)-f^\star)^2$ for $k\ge k_0$, where $k_0$ is the first $k\ge b$ with $f(x_k)-f^\star<\frac12\mu\omega^2$.
\end{theorem}
\begin{proof}
By Lemma~\ref{lem:descent}, $f(x_k)$ decreases to $f^\star$, so $k_0$ exists and $f(x_k)-f^\star<\frac12\mu\omega^2$ for $k\ge k_0$. Fix such a $k$ and shift the indices so that $k=0$. If some $\nabla f(x_i)$, $0\le i\le m$, vanishes, the claim is trivial, so let them be nonzero. Strong convexity and $f(x_i)\le f(x_0)$ give $\norm{x_i-x^\star}^2\le2(f(x_0)-f^\star)/\mu<\omega^2$, so all segments between $x_0,\dots,x_m,x^\star$ lie in $B$.

By Lemma~\ref{lem:descent}, $\nabla f(x_i)\perp \nabla f(x_{i-1}),s_i$, so strong convexity gives $\mu\norm{s_i}^2\le\ip{\nabla f(x_i)-\nabla f(x_{i-1})}{s_i}\le\norm{\nabla f(x_{i-1})}\norm{s_i}$. Put $z_i=\nabla f(x_i)-\nabla f(x_{i-1})-\chi s_i$ for $1\le i\le m$. Then $z_i=Rs_i+\epsilon_i$ with
\[
\epsilon_i=\int_0^1\bigl(\nabla^2f(x_{i-1}+ts_i)-\nabla^2f(x^\star)\bigr)s_i\,dt,\qquad\norm{\epsilon_i}^2\le\frac{2M^2}\mu(f(x_0)-f^\star)\norm{s_i}^2.
\]
For $1\le i<j\le m$ the space $S_{j-1}$ covers $\nabla f(x_{i-1})$, $\nabla f(x_i)$ and $s_i$, so $\nabla f(x_j)\perp z_i$, and $\nabla f(x_1),\dots,\nabla f(x_m)$ are mutually orthogonal. Let $Q$ have the columns $\nabla f(x_j)/\norm{\nabla f(x_j)}$ and $\mathsf X$ the columns $Q^\top z_j/\norm{\nabla f(x_{j-1})}$, $1\le j\le m$. This matrix is upper triangular with the diagonal entries $\norm{\nabla f(x_j)}/\norm{\nabla f(x_{j-1})}$, so $\det\mathsf X=\norm{\nabla f(x_m)}/\norm{\nabla f(x_0)}$. Its columns are $a_j+e_j$ with $a_j=Q^\top Rs_j/\norm{\nabla f(x_{j-1})}$ and $e_j=Q^\top\epsilon_j/\norm{\nabla f(x_{j-1})}$. Write $\mathsf X=A+E$, where $A$ and $E$ have columns $a_j$ and $e_j$. Since $\operatorname{rank}A\le\operatorname{rank}R<m$, we have $\det A=0$. The bounds $\norm{a_j}\le\norm R/\mu$ and $\norm{e_j}\le(2M^2(f(x_0)-f^\star)/\mu^3)^{1/2}\le M\omega/\mu$ place both matrices in a bounded set of $m\times m$ matrices. The determinant is Lipschitz on that set, so
\[
\frac{\norm{\nabla f(x_m)}}{\norm{\nabla f(x_0)}}
=|\det(A+E)-\det A|
\le C\norm E\le C'\sqrt{f(x_0)-f^\star},
\]
with constants depending only on $m$, $\mu$, $M$, $\norm R$ and $\omega$. Squaring and using $f(x_m)-f^\star\le\norm{\nabla f(x_m)}^2/(2\mu)$ and $\norm{\nabla f(x_0)}^2\le2L(f(x_0)-f^\star)$ proves the claim.
\end{proof}

On a quadratic with $\operatorname{rank}(\nabla^2f-\chi I)\le m-1$, $M=0$ and every method of the theorem terminates within $m$ steps. The constant is not optimized, so the theorem concerns the order of convergence rather than the size of the quadratic regime. Near the minimizer of the instances of Proposition~\ref{prop:smooth}, $\nabla^2G=I+(K-1)ee^\top$, and covering with $m=2$ needs only $\nabla f(x_k)$, $\nabla f(x_{k-1})$ and $s_k$. SESOP with the previous gradient added therefore reaches the minimizer two steps after entering a ball in the quadratic region, while SESOP itself keeps $G(x_{k+1})\ge(K-1)G(x_k)/(K+127)$, so adding one vector turns the rate of gradient descent into finite termination.

Once a window of length $m$ is filled, its $m$ most recent gradients are mutually orthogonal, so $\dim S_k\ge m$ before termination. For the pure window method also $\dim S_k\le2m$. Direct formation of the reduced Hessian uses one Hessian-vector product per basis vector, hence $\Theta(\sqrt\kappa)$ such products per inner Newton iteration when $m\asymp\sqrt\kappa$. Our bounds count exact subspace steps and do not bound the number of inner Newton iterations. Section~\ref{sec:filters} stores a horizon of this order in two vectors.

\subsection{Filters on the right time scale}\label{sec:filters}

Two geometric filters keep the virtual momentum of Lemma~\ref{lem:cert} for a constant step. Any $a\le1/\sqrt\kappa$ is admissible there, with the rate $1-a$, but the best choice $a=1/\sqrt\kappa$ is unknown. With $\pi_b=0$ and $a_k=a$, unrolling \eqref{eq:prec} gives
\[
\pi_k=y_k+\frac a\lambda z_k,\qquad z_k=\sum_{i=b+1}^{k}\rho^{k-i}\nabla f(x_i),\qquad y_k=\sum_{i=b+1}^{k}\rho^{k-i+1}s_i ,
\]
with the factor $\rho=1-a$. A method that stores $z_k$ and $y_k$ separately never needs the coefficient $a/\lambda$, which involves the unknown constants, and has to choose only $\rho$.

We restrict $a$ to the dyadic values $\alpha_J=2^{-J}$ and take $\lambda=4^{-J}L$. A level $J\ge1$ is admissible exactly when $4^J\ge\kappa$, and it is nearly optimal when also $4^J\le4\kappa$. With $\rho_J=1-\alpha_J$, the \emph{filter pair} of level $J$ started at $b$ is
\begin{equation}\label{eq:pair}
z_b=y_b=0,\qquad z_{k+1}=\rho_Jz_k+\nabla f(x_{k+1}),\qquad y_{k+1}=\rho_J\bigl(y_k+s_{k+1}\bigr)\qquad(k\ge b).
\end{equation}

\begin{lemma}[certificate of a filter pair]\label{lem:pair}
Let $f\in\Fc_{\mu,L}$ and $\kappa\le4^J$, and let the iterates satisfy \eqref{eq:step} with any rule for $S_k$. If the pair \eqref{eq:pair} of level $J$ started at $b$ lies in $S_k$ for $b\le k<b+h$, then $f(x_{b+h+1})-f^\star\le2\rho_J^h(f(x_b)-f^\star)$.
\end{lemma}
\begin{proof}
With $\alpha=2^{-J}$ and $\lambda=\alpha^2L\le\mu$ we have $\alpha^2/\lambda=1/L$. Substitution in \eqref{eq:pair} shows that $\pi_k=y_k+\frac\alpha\lambda z_k$ satisfies \eqref{eq:prec} with $a_k=\alpha$ and $\pi_b=0$. Since $\pi_k\in S_k$, Lemma~\ref{lem:cert} with $C_b\le2(f(x_b)-f^\star)$ gives the claim.
\end{proof}

For a matching level $\kappa\le4^J\le4\kappa$, the pair has $\delta=2^{-J}\ge1/(2\sqrt\kappa)$, as Corollary~\ref{cor:memA} requires for a bounded window, and the rate is $\exp(-1/(2\sqrt\kappa))$ per step with two stored vectors, where the window of Proposition~\ref{prop:window} needs about $2.8\sqrt\kappa$, that is thousands of vectors for $\kappa=10^6$.

\paragraph{Activation schedule.} Algorithm~\ref{alg:schedule} creates the pair of level $j$ at iteration $2^j$, when the run is as long as its horizon $1/\alpha_j$, and never discards it. It uses no constant of the function and no test.

\begin{algorithm}[t]
\caption{Activation schedule}\label{alg:schedule}
\begin{algorithmic}[1]
\Require starting point $x_0$
\State Compute $\nabla f(x_0)$.
\For{$k=0,1,2,\dots$}
\State If $k=2^j$ with $j\ge1$, create the pair $z^{(j)}_k=y^{(j)}_k=0$.
\State Set $S_k=\spn\bigl(\{\nabla f(x_k)\}\cup\{z^{(j)}_k,y^{(j)}_k:2\le2^j\le k\}\bigr)$.
\State Compute $x_{k+1}=\argmin_{x_k+S_k}f$, $s_{k+1}=x_{k+1}-x_k$ and $\nabla f(x_{k+1})$.
\State Update every pair by \eqref{eq:pair} with its own factor $\rho_j$.
\EndFor
\end{algorithmic}
\end{algorithm}

\begin{theorem}[activation schedule]\label{thm:schedule}
Let $f\in\Fc_{\mu,L}$. Algorithm~\ref{alg:schedule} satisfies $f(x_{k+1})-f^\star\le2e\exp(-k/(2\sqrt\kappa))(f(x_0)-f^\star)$ for $k\ge0$, and at iteration $k\ge2$ it stores $2\lfloor\log_2k\rfloor$ vectors. Hence $N(\varepsilon)\le1+\lceil2\sqrt\kappa\log(2e/\varepsilon)\rceil$, with $O(1+\log\kappa+\log(1+\log(1/\varepsilon)))$ vectors stored by then.
\end{theorem}
\begin{proof}
Let $J=\max\{1,\lceil\log_4\kappa\rceil\}$, so that $\kappa\le4^J\le4\kappa$ and $\alpha=2^{-J}$ satisfies $\frac1{2\sqrt\kappa}\le\alpha\le\frac1{\sqrt\kappa}$. The pair of level $J$ is activated at $\tau=2^J=1/\alpha\le2\sqrt\kappa$ with zero state. Lemma~\ref{lem:pair} with origin $b=\tau$ and monotonicity give, for $k\ge\tau$,
\[
\begin{aligned}
f(x_{k+1})-f^\star
&\le2(1-\alpha)^{k-\tau}(f(x_\tau)-f^\star)\\
&\le2e^{1-\alpha k}(f(x_0)-f^\star)\\
&\le2e\,e^{-k/(2\sqrt\kappa)}(f(x_0)-f^\star).
\end{aligned}
\]
For $k<\tau$ the right-hand side exceeds $2(f(x_0)-f^\star)$, and the bound holds by monotonicity. At iteration $k\ge2$ the active pairs are those with $2\le2^j\le k$, which gives the storage, and solving the rate for $k$ gives $N(\varepsilon)$.
\end{proof}

The schedule thus attains the accelerated upper bound $O(\sqrt\kappa\log(1/\varepsilon))$ in exact subspace steps without problem constants. At iteration $k\ge2$ its subspace has dimension at most $1+2\lfloor\log_2k\rfloor$. The bound counts each subspace minimization as one step; it does not bound the number of elementary oracle calls needed by an inner solver. Inexact solves and the per-step cost in the numerical examples are addressed in Sections~\ref{sec:inexact} and~\ref{sec:experiments}. Not knowing $\kappa$ costs this construction a factor two in the exponent and memory that grows like $\log k$, with no restart, test or estimate of a constant involved.

\paragraph{Smooth adaptive pair.} Algorithm~\ref{alg:smooth} updates one pair with the factor $\rho_k=1-1/H_k$. The scalar $R_k$ averages $\norm{\nabla f(x_k)}^2$ with the factor $1-1/(2H_k)$, and $\ell_k=\log(R_{k+1}/((1-\frac1{2H_k})R_k))$ measures how much the gradient energy fails to decay on the current horizon. The horizon never decreases, and $(\log H_k)^2$ grows by $\ell_k/4$, so $H_k$ grows quickly only while this happens.

\begin{algorithm}[t]
\caption{Smooth adaptive pair}\label{alg:smooth}
\begin{algorithmic}[1]
\Require starting point $x_0$
\State Compute $\nabla f(x_0)$; if it vanishes, \Return $x_0$.
\State Set $H_0=2$, $R_0=\norm{\nabla f(x_0)}^2$ and $z_0=y_0=0$.
\For{$k=0,1,2,\dots$}
\State Set $a_k=1/H_k$, $\rho_k=1-a_k$ and $S_k=\spn\{\nabla f(x_k),z_k,y_k\}$.
\Statex\hspace{\algorithmicindent}Optionally include further directions in $S_k$.
\State Compute $x_{k+1}=\argmin_{x_k+S_k}f$, $s_{k+1}=x_{k+1}-x_k$ and $\nabla f(x_{k+1})$.
\State Update $z_{k+1}=\rho_kz_k+a_k\nabla f(x_{k+1})$ and $y_{k+1}=\rho_k(y_k+s_{k+1})$.
\State Update $R_{k+1}=(1-\frac{a_k}2)R_k+\frac{a_k}2\norm{\nabla f(x_{k+1})}^2$.
\Statex\hspace{\algorithmicindent}Set $\ell_k=\log\bigl(1+\frac{a_k\norm{\nabla f(x_{k+1})}^2}{(2-a_k)R_k}\bigr)$.
\State Update $H_{k+1}=\exp\bigl(\sqrt{(\log H_k)^2+\ell_k/4}\bigr)$.
\EndFor
\end{algorithmic}
\end{algorithm}

\begin{theorem}[smooth adaptive pair]\label{thm:afb}
Let $f\in\Fc_{\mu,L}$ with $f(x_0)-f^\star>0$, and let the iterates of Algorithm~\ref{alg:smooth} satisfy \eqref{eq:step}, with or without further directions in $S_k$. The adaptive pair requires two stored vectors and two scalars, which are never reset; additional directions may require further storage. There are absolute constants $C\ge1$ and $c>0$ such that
\[
\begin{aligned}
f(x_k)-f^\star
&\le C\kappa^5e^{-ck/\sqrt\kappa}(f(x_0)-f^\star),\\
N(\varepsilon)
&\le\Bigl\lceil c^{-1}\sqrt\kappa\Bigl[\log\frac1\varepsilon+5\log\kappa+\log C\Bigr]\Bigr\rceil.
\end{aligned}
\]
If $S_k$ is spanned by the pair, $\nabla f(x_k)$, $s_k$, the anchor $x_k-x_0$ and Nemirovski's aggregate with $w_0=1$, Proposition~\ref{prop:sesop} applies as well.
\end{theorem}
\begin{proof}
The proof has three steps. The stored vectors are bounded by the gradient energy $R_k$, the horizon is bounded by $C_H\sqrt\kappa$ for an absolute constant $C_H\ge1$ when $\kappa\ge4$, and a summation by parts turns the growth of $H_k$ into the rate. Put $q_k=1-a_k/2$. Since $f(x_0)>f^\star$, we have $R_0>0$. Starting from $H_0=2$, the updates give inductively $H_{k+1}\ge H_k\ge2$ and $R_{k+1}\ge\frac34R_k>0$. Thus all logarithms below are well defined, including after a minimizer is reached. By the update of $R_k$, $R_{k+1}=q_kR_ke^{\ell_k}$, and the update of $H_k$ gives
\begin{equation}\label{eq:H-identity}
4(\log H_{k+1})^2-4(\log H_k)^2=\ell_k\ge0 .
\end{equation}
Hence $H_k$ is nondecreasing, $H_k\ge2$, $a_k\le\frac12$ and $q_k\ge\frac34$.

\emph{Step 1: the stored states.} We show that for every $k\ge0$
\begin{equation}\label{eq:state-bounds}
\begin{aligned}
\norm{z_k}^2&\le R_k,\qquad
\norm{y_k}^2\le\tfrac{2H_k^2}{\mu^2}R_k,\\
\norm{\nabla f(x_k)}^2&\le2H_kR_k,\\
2\mu(f(x_k)-f^\star)&\le R_k\le2L(f(x_0)-f^\star).
\end{aligned}
\end{equation}
Since $z_k\in S_k$ and $\nabla f(x_{k+1})\perp S_k$, $\norm{z_{k+1}}^2=\rho_k^2\norm{z_k}^2+a_k^2\norm{\nabla f(x_{k+1})}^2\le R_{k+1}$ by induction, because $\rho_k^2\le q_k$ and $a_k^2\le a_k/2$. For $y_k$ put $c_j=\prod_{t=j}^{k-1}\rho_t$, so that $y_k=\sum_{j=0}^{k-1}c_js_{j+1}$. Strong convexity and $\nabla f(x_{j+1})\perp s_{j+1}$ give $\mu\norm{s_{j+1}}^2\le\ip{\nabla f(x_{j+1})-\nabla f(x_j)}{s_{j+1}}\le\norm{\nabla f(x_j)}\norm{s_{j+1}}$. Since $H_t\le H_k$ for $t\le k$, $c_j\le(1-1/H_k)^{k-j}$ and $\sum_jc_j\le H_k$. Expanding the recursion of $R_k$, the coefficient of $\norm{\nabla f(x_0)}^2$ is $\prod_{t<k}q_t\ge c_0$, and that of $\norm{\nabla f(x_j)}^2$, $1\le j<k$, is $\frac{a_{j-1}}2\prod_{t=j}^{k-1}q_t\ge c_j/(2H_k)$, because $q_t\ge\rho_t$ and $a_{j-1}\ge1/H_k$. Hence $2H_kR_k\ge\sum_{j<k}c_j\norm{\nabla f(x_j)}^2$, and the Cauchy--Schwarz inequality gives
\[
\norm{y_k}\le\frac1\mu\sum_{j=0}^{k-1}c_j\norm{\nabla f(x_j)}\le\frac1\mu\Bigl(\sum_jc_j\Bigr)^{1/2}\Bigl(\sum_jc_j\norm{\nabla f(x_j)}^2\Bigr)^{1/2}\le\frac{\sqrt2H_k}\mu\sqrt{R_k}.
\]
For $k\ge1$, $R_k\ge\frac{a_{k-1}}2\norm{\nabla f(x_k)}^2\ge\norm{\nabla f(x_k)}^2/(2H_k)$, and $R_0=\norm{\nabla f(x_0)}^2$, which gives the third bound. Finally $R_k$ is a convex combination of $\norm{\nabla f(x_j)}^2$, $j\le k$, and $2\mu(f(x_k)-f^\star)\le\norm{\nabla f(x_j)}^2\le2L(f(x_0)-f^\star)$ for $j\le k$ by Lemma~\ref{lem:descent}.

\emph{Step 2: the horizon.} All implicit constants in this step are absolute. Let $\kappa\ge4$. Lemma~\ref{lem:descent} and \eqref{eq:state-bounds} give $\norm{\nabla f(x_{k+1})}^2\le\kappa\norm{\nabla f(x_k)}^2\le2\kappa H_kR_k$, hence $\ell_k\le\log(1+\frac43\kappa)$. If $H_k<\sqrt\kappa$ for all $k$, the desired bound holds. Otherwise let $b$ be the first index with $H_b\ge\sqrt\kappa$. By \eqref{eq:H-identity},
\[
4(\log H_b)^2\le(\log\kappa)^2+\log(1+\tfrac43\kappa)
=(\log\kappa)^2+\log\kappa+O(1),
\]
so $H_b=O(\sqrt\kappa)$.

For $t\ge b$ we have $a_t\le1/\sqrt\kappa$. The vector $\pi_t=y_t+z_t/\mu$ lies in $S_t$ and satisfies \eqref{eq:prec}, so Lemma~\ref{lem:cert} applies with $\lambda=\mu$. By \eqref{eq:state-bounds} and strong convexity, its initial constant satisfies $C_b=O(H_b^2R_b/\mu)$. Thus $\norm{\nabla f(x_{t+1})}^2\le C\kappa H_b^2R_b\prod_{i=b}^{t-1}\rho_i$ for an absolute $C$. Put $P_t=\exp(\sum_{i=b}^{t-1}\ell_i)$ and $V_t=\prod_{i=b}^{t-1}(\rho_i/q_i)$. Since $R_t=R_bP_t\prod_{i=b}^{t-1}q_i$ and $1-\rho_t/q_t=a_t/(2q_t)$,
\[
P_{t+1}-P_t
=\frac{a_t\norm{\nabla f(x_{t+1})}^2}{2q_tR_b\prod_{i=b}^{t-1}q_i}
\le C\kappa H_b^2(V_t-V_{t+1}).
\]
Telescoping gives $\sum_{i=b}^{t-1}\ell_i\le\log(1+C\kappa H_b^2)\le2\log\kappa+O(1)$. Therefore \eqref{eq:H-identity} yields $4(\log H_t)^2\le(\log\kappa)^2+3\log\kappa+O(1)$, and hence $H_t\le C_H\sqrt\kappa$ uniformly for $\kappa\ge4$.

\emph{Step 3: the rate.} If $\kappa<4$, Lemma~\ref{lem:descent} gives $f(x_k)-f^\star\le e^{-k/\kappa}(f(x_0)-f^\star)\le e^{-k/(2\sqrt\kappa)}(f(x_0)-f^\star)$. Let $\kappa\ge4$ and $k\ge1$. Then $\log(R_{i+1}/R_i)=\log q_i+\ell_i$ and $H_i\log q_i\le-\frac12$. By \eqref{eq:H-identity} and $H_i\le t$ on $[H_i,H_{i+1}]$,
\[
\sum_{i=0}^{k-1}H_i\ell_i=\sum_{i=0}^{k-1}\int_{H_i}^{H_{i+1}}\frac{8H_i\log t}t\,dt\le\int_{H_0}^{H_k}8\log t\,dt\le8H_k\log H_k .
\]
Put $M_0=2L(f(x_0)-f^\star)$ and $u_i=\log(R_i/M_0)\le0$. Summation by parts gives
\[
\sum_{i=0}^{k-1}H_i(u_{i+1}-u_i)
=H_ku_k-H_0u_0-\sum_{i=1}^{k}(H_i-H_{i-1})u_i\ge H_ku_k,
\]
since $H_i\ge H_{i-1}$ and $u_i\le0$. Combining these bounds with $\log(R_{i+1}/R_i)=\log q_i+\ell_i$, exponentiating and using $R_k\ge2\mu(f(x_k)-f^\star)$ gives $(f(x_k)-f^\star)/(f(x_0)-f^\star)\le\kappa H_k^8e^{-k/(2H_k)}$, and $H_k\le C_H\sqrt\kappa$ gives the claimed bound with $C=C_H^8$ and $c=1/(2C_H)$. The case $k=0$ is immediate since $C\kappa^5\ge1$. Solving the bound for $k$ gives $N(\varepsilon)$. Only the pair enters the proof, so further directions are allowed, and as the pair consists of combinations of past gradients and steps, Proposition~\ref{prop:sesop} applies to the SESOP variant.
\end{proof}

The prefactor is polynomial in $\kappa$, so \eqref{eq:accel} is not established, and we do not know if it can be removed. Unlike $\mu$-estimates in restart schemes \cite{nesterov2013gradient,lin2015adaptive}, $H_k$ only sets the forgetting factor; Step~2 bounds it by an absolute multiple of $\sqrt\kappa$ for $\kappa\ge4$.

\begin{proof}[of Remark~\ref{rem:attenuation}, periodic attenuation]
Put $a=1/h$, $\lambda=L/h^2\le\mu$ and $W_b=f(x_b)-f^\star+\lambda\norm{y_b}^2+L^{-1}\norm{z_b}^2$ at the start $b$ of a block. The virtual momentum $\pi_k=y_k+\frac hLz_k$ lies in $S_k$ and satisfies \eqref{eq:prec}. Strong convexity and Lemma~\ref{lem:cert} therefore give
\[
f(x_{b+T})-f^\star\le C_1\rho^{T-1}W_b\le C_1e^{-A}W_b.
\]
Before attenuation, $x_{b+T}-y_{b+T}$ is a convex combination of $x_b-y_b$ and the iterates within the block, because $x_{k+1}-y_{k+1}=(1-\rho)x_{k+1}+\rho(x_k-y_k)$. The iterates within the block lie within $\sqrt{2(f(x_b)-f^\star)/\mu}$ of $x^\star$. Combining this bound with the initial term $y_b$ and the recursion for $z_k$, where $\sum_{j\ge0}\rho^j=h$, gives
\[
\lambda\norm{y_{b+T}}^2+L^{-1}\norm{z_{b+T}}^2\le C_2h^2W_b
\]
for absolute constants $C_1,C_2$ independent of the block length. After multiplication of both registers by $\eta/h$, the potential is at most $(C_1e^{-A}+C_2\eta^2)W_b$. Choose $A$ sufficiently large and then $\eta\in(0,1)$ sufficiently small to make this factor at most $1/2$. Since $W_0=f(x_0)-f^\star$, monotonicity and $T=O(\sqrt\kappa)$ give the claimed accelerated rate. Normalized weights grow by $\rho^{-1}$ per step and $(h/\eta)$ per block, so $\delta=1-\rho(\eta/h)^{1/T}\asymp\log\kappa/\sqrt\kappa$.
\end{proof}

\subsection{Inexact subproblems}\label{sec:inexact}

Let the gradients and the filter updates be exact, fix $\gamma\in[0,1]$, put $a_k=2^{-\lfloor\log_2k\rfloor}$ for $k\ge2$ and $a_0=a_1=\frac12$, write $\nu_k=\norm{P_{S_k}\nabla f(x_{k+1})}$ for the reduced gradient at the computed point, and replace \eqref{eq:step} by
\begin{equation}\label{eq:inexact}
x_{k+1}\in x_k+S_k,\qquad f(x_k)-f(x_{k+1})\ge\frac{\norm{\nabla f(x_k)}^2}{2L'},\qquad\nu_k^2\le\gamma a_k\norm{\nabla f(x_k)}^2,
\end{equation}
with $\nabla f(x_k)\in S_k$ and a constant $L'\ge L$ that enters only the analysis, and put $\kappa'=L'/\mu$. A monotone inner solver started at a gradient step with sufficient decrease satisfies the second condition, and the third is observable.

Paying the defect $E_k$ of Lemma~\ref{lem:cert} from the decrease of $f$ slows a pair only slightly.
\begin{lemma}[decrease budget]\label{lem:budget}
In the setting of Lemma~\ref{lem:cert} let $a_k=a\in(0,1)$, $\rho=1-a$, and $E_k\le q(f(x_k)-f(x_{k+1}))$ for $k\ge b$ with $q\ge0$. Put $r_q=1-a/(1+aq)$ and $\beta_q=a\rho q/r_q$. Then $\beta_q\le aq$ and $f(x_{b+h})-f^\star\le r_q^{h-1}(C_b+\beta_q(f(x_b)-f^\star))$ for $h\ge1$.
\end{lemma}
\begin{proof}
Put $t=a\rho q$ and $U_j=\rho^jC_b+e_{b+j}$. Lemma~\ref{lem:cert} gives $f(x_{b+j+1})-f^\star\le U_j$, and the recursion of $e_k$ gives
\[
U_{j+1}=\rho U_j+a\rho E_{b+j}\le\rho U_j+t\bigl(f(x_{b+j})-f(x_{b+j+1})\bigr).
\]
Since $\rho\le r_q\le1$, we have $t=r_q\beta_q\le\beta_q\le t/\rho=aq$, and $\beta_q-t=\beta_q(1-r_q)\le a^2q/(1+aq)=r_q-\rho$. For $V_j=U_j+\beta_q(f(x_{b+j})-f^\star)$ these relations and $f(x_{b+j+1})-f^\star\le U_j$ give
\begin{align*}
V_{j+1}&\le\rho U_j+t(f(x_{b+j})-f^\star)+(\beta_q-t)(f(x_{b+j+1})-f^\star)\\
&\le(\rho+\beta_q-t)U_j+t(f(x_{b+j})-f^\star)\le r_qV_j .
\end{align*}
With $f(x_{b+h})-f^\star\le U_{h-1}\le V_{h-1}$ and $V_0=C_b+\beta_q(f(x_b)-f^\star)$ this gives the bound.
\end{proof}

\begin{proposition}[inexact activation schedule]\label{thm:inexact}
Let $f\in\Fc_{\mu,L}$ with $f(x_0)-f^\star>0$, and let the iterates of Algorithm~\ref{alg:schedule} satisfy \eqref{eq:inexact}. Then, for $k\ge0$,
\[
f(x_{k+1})-f^\star
\le(2+\gamma)e^{1/(1+\gamma)}
\exp\Bigl(-\frac{k}{2(1+\gamma)\sqrt{\kappa'}}\Bigr)
(f(x_0)-f^\star),
\]
also with additional directions in $S_k$.
\end{proposition}
\begin{proof}
Let a pair of level $J$ with $a=2^{-J}$ and $\lambda=a^2L'\le\mu$ lie in $S_k$. As in Lemma~\ref{lem:pair}, its virtual momentum $\pi_k$ lies in $S_k$, and $s_{k+1}\in S_k$ by \eqref{eq:inexact}, so $\pi_k+s_{k+1}\in S_k$. Lemma~\ref{lem:cert} and \eqref{eq:inexact} give
\[
E_k\le\frac{\nu_k^2}{2\lambda}\le\frac{\gamma a_k\norm{\nabla f(x_k)}^2}{2a^2L'}\le\frac{\gamma a_k}{a^2}\bigl(f(x_k)-f(x_{k+1})\bigr).
\] If $a_k\le a$, Lemma~\ref{lem:budget} applies with $q=\gamma/a$, so that $aq=\gamma$. A pair created at $b$ has $\pi_b=0$ and $C_b\le2(f(x_b)-f^\star)$, hence
\begin{equation}\label{eq:inexact-pair}
f(x_{b+h})-f^\star\le(2+\gamma)\Bigl(1-\frac a{1+\gamma}\Bigr)^{h-1}(f(x_b)-f^\star)\qquad(h\ge1).
\end{equation}
Choose $J=\max\{1,\lceil\log_4\kappa'\rceil\}$ as in Theorem~\ref{thm:schedule}, so that $a\ge1/(2\sqrt{\kappa'})$ and $\lambda\le\mu$. Its activation time $b=1/a$ ensures $a_k\le a$ for $k\ge b$. The same activation-time and monotonicity argument, applied to \eqref{eq:inexact-pair} with decay $a/(1+\gamma)$ and prefactor $2+\gamma$, gives the stated bound.
\end{proof}

Unlike \cite{kuruzov2024some}, where the gradients are inexact, here the defect of Lemma~\ref{lem:cert} is charged to the actual decrease, so the stopping rule is observable: at the matching level the relative residual may be of order $\kappa'^{-1/4}$, and the tolerance decreases like $k^{-1/2}$. Sufficient decrease alone certifies only the rate $1-1/\kappa'$ of the gradient step. For $\gamma=0$ and $L'=L$ the proposition reduces to Theorem~\ref{thm:schedule}.

\section{Numerical experiments}\label{sec:experiments}

The experiments address two questions. Does SESOP, run as a black box, follow the slow trajectories of Section~\ref{sec:main}, and what do the guarantees of Section~\ref{sec:upper} cost outside the worst case? All experiments use Python~3.13 with NumPy~2.5 and SciPy~1.18. Times are single runs with one BLAS thread, so differences of a few percent are insignificant.

\begin{figure}[t]
\centering
\includegraphics[width=0.88\textwidth]{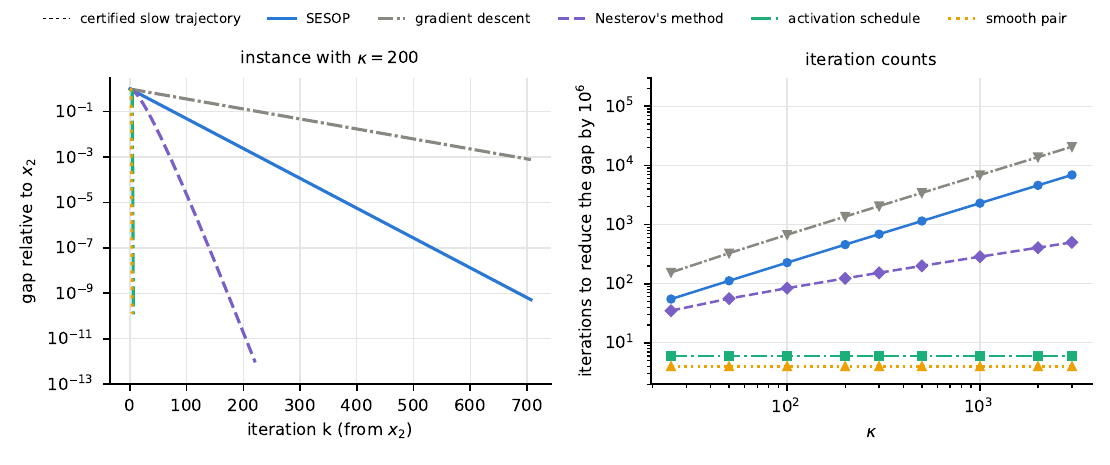}
\caption{Black-box runs on the SESOP instances from $x_2$. Left: $\kappa=200$, SESOP follows the prescribed slow trajectory (dashed). Right: iterations to reduce the gap by $10^6$ against $\kappa$. Slopes fit $\text{iterations}\approx c\kappa^s$ by least squares. Gradient descent uses the step $1/L$, Nesterov's method the true $\mu$ and $L$.}\label{fig:hard}
\end{figure}

\subsection{Hard instances}

For SESOP with one aggregate, the linear weights and $\kappa\in[25,3000]$ we computed a five-dimensional variant of the construction, with two explicit initial points instead of the orthogonal segment, and the methods evaluate its interpolant (Lemma~\ref{lem:interp}) by a quadratic program over the simplex. The conditions of the proof of Lemma~\ref{lem:bridge} are solved by backward substitution for $h_n$, $2\le n\le N$, with $N$ from $150$ to $8000$, beyond every SESOP run, and $\vartheta$ from $10^{-9}$ to $10^{-3}$, with residual below $10^{-14}$. All interpolation residuals are positive, the smallest relative one between $6\cdot10^{-9}$ and $1.3\cdot10^{-2}$. The data are not certified in exact arithmetic, so the runs only illustrate Theorem~\ref{thm:slow}.

For $\kappa\ge50$ double precision cannot resolve the first line search from $f_0=12\vartheta^{-2}>10^8$, so the runs in Figure~\ref{fig:hard} start at $x_2$, SESOP with its state after two iterations and the filter methods with empty memory. SESOP follows the prescribed trajectory with relative deviation below $4\cdot10^{-9}$ and the rate $-2\log(1-3/\kappa)$ to thirteen digits. Its counts grow like $\kappa^{1.01}$, those of gradient descent like $\kappa^{1.02}$ and those of Nesterov's method like $\kappa^{0.55}$. The activation schedule and the smooth pair leave the prescribed path at once and need $6$ and $4$ iterations, which says nothing about their worst case.

Lemma~\ref{lem:invariant} predicts the factor $r^2$ per step. The plane rule of Proposition~\ref{prop:plane3} keeps this factor to twelve digits for $3000$ iterations. For the smooth SESOP instances of Proposition~\ref{prop:smooth}, we check subspace membership and gradient orthogonality for the three prescribed initial steps, then run SESOP on the quadratic tail from $x_3$. With Nemirovski's and the linear weights, $K\in\{30,100,300,1000\}$ and $\theta=10^{-8}$, the tail keeps the factor to six digits over a run of $1500$ iterations including the three initial steps. These checks confirm the prescribed slow behavior in floating point.

\subsection{Logistic regression and capillary energies}

We use regularized logistic regression $f(x)=\frac1\ell\sum_{i=1}^{\ell}\log(1+e^{-y_i\ip{a_i}{x}})+\frac\mu2\norm{x}^2$ with $\ell=2000$, $n=1000$, a data matrix obtained from a Gaussian matrix by replacing its singular values with logarithmically spaced values spanning two decades, labels from a noisy linear model, and $L=1$. A harder test is the capillary-surface problem \cite{finn1986equilibrium}: a liquid surface $u$ over $\Omega=(0,1)^d$ minimizes its area plus the energy of gravity and of a load $q$, $E(u)=\int_\Omega(\sqrt{1+\norm{\nabla u}^2}+\frac12u^2-qu)\,dx$ over $u\in H^1_0(\Omega)$. The gravity term makes $E$ strongly convex in the $L^2(\Omega)$ norm, and the curvature $(1+t^2)^{-3/2}$ of the area integrand degenerates where the surface is steep, so under a strong load the Hessian varies by orders of magnitude. For $d=1$ we use finite differences on $n$ interior points, so that $\kappa\le4(n+1)^2+1$, and $q(x)=Q(\sin3\pi x+\frac12\operatorname{sign}\sin7\pi x+\frac3{10}\mathbf 1_{\{x>3/10\}})$. All methods start at zero, and $f^\star$ is computed by Newton's method.

Nesterov's method (AGD) uses $\mu$ and $L$, its restarted versions \cite{odonoghue2015adaptive} needed between $0.65$ and $1.9$ times its iterations, and L-BFGS is the SciPy code with memory $10$ \cite{liu1989limited,nocedal2006numerical}. SESOP uses one aggregate with Nemirovski's weights. The experimental window uses the $8$ most recent pairs $(\nabla f(x_i),s_{i+1})$ and the current gradient, giving nine gradients and eight steps once filled; it thus includes one more gradient than a window of length $m=8$ in Section~\ref{sec:prelim}. The subspace problems are solved by damped Newton iterations, with reduced-gradient stopping tolerances $10^{-10}\norm{\nabla f(x_k)}$ for logistic regression and $10^{-12}\norm{\nabla f(x_k)}$ for the one-dimensional capillary problems.

\begin{table}[t]
\caption{Iterations to reduce $f(x_k)-f^\star$ by $10^{10}$ (logistic) and $10^8$ (capillary, $d=1$, $n=400$).}\label{tab:scaling}
\footnotesize\setlength{\tabcolsep}{3pt}
\begin{tabular*}{\textwidth}{@{\extracolsep\fill}lrrrrrr}
\toprule
& \multicolumn{3}{c}{logistic, $\kappa$} & \multicolumn{3}{c}{capillary, $Q$}\\
\cmidrule{2-4}\cmidrule{5-7}
method & $10^2$ & $10^4$ & $10^6$ & $1$ & $10$ & $40$\\
\midrule
AGD (true $\mu$) & 101 & 1155 & 12469 & 6831 & 9332 & 9415\\
L-BFGS & 38 & 131 & 179 & 1078 & 4441 & 4463\\
plane CG & 36 & 123 & 161 & 556 & 4351 & 6091\\
SESOP & \textbf{35} & 118 & 167 & 450 & 3259 & 4724\\
window, $8$ pairs & \textbf{35} & \textbf{114} & \textbf{159} & 445 & 3755 & 3930\\
activation schedule & \textbf{35} & 116 & 162 & \textbf{435} & 3188 & \textbf{3147}\\
smooth pair & 36 & 121 & 174 & 933 & 4427 & 4554\\
smooth pair + SESOP dir. & \textbf{35} & 117 & 165 & 448 & \textbf{3183} & 3628\\
\bottomrule
\end{tabular*}
\end{table}

On logistic regression all subspace methods behave like accelerated methods, and their counts grow much more slowly than $\sqrt\kappa$, since the condition number near the solution is far below the nominal one (Table~\ref{tab:scaling}). The methods with guarantees need at most $10$ percent more iterations than the best one. The horizon of the smooth pair stays below $5$ there, and on the capillary energy with $n=200$ below $0.09\sqrt\kappa$, with $\kappa$ of the Hessian at the solution. For the capillary results with $n=400$ in Table~\ref{tab:scaling}, the activation schedule is the fastest method or within $0.2$ percent of it for every load, with a search space of dimension about $20$ against $4$ for SESOP, and the SESOP directions remove most of the loss of the smooth pair, which alone needs up to $2.3$ times more iterations. On these problems the guarantee of the schedule thus costs almost nothing.

\subsection{Cost of an iteration and inexact subproblems}

For $f(x)=\sum_{i=1}^{\ell}\varphi_i(\ip{a_i}{x})+\frac\mu2\norm x^2$ the reduced Hessian needs only the images $AD_k$ of the stored directions, and if they follow the same recursions, a step costs one product with $A$ and one with $A^\top$ \cite{narkiss2005sequential,zibulevsky2010l1}. We timed logistic regression with $\ell=10000$, $n=2000$, $\kappa=10^5$ and inner tolerance $10^{-4}\norm{\nabla f(x_k)}$, and the capillary energy on $(0,1)^2$ with $Q=40$ in $64^2$ Galerkin sine modes, that is $\Phi(Ac)$ with a dense $19522\times4096$ matrix $A$ and $\kappa\approx2.4\cdot10^4$. To reduce the gap by $10^{10}$ on logistic regression the subspace methods needed $264$ to $280$ iterations and $6.4$ to $7.9$ seconds, L-BFGS $305$ iterations and $10.1$ seconds, and AGD with or without gradient restart about $3700$ iterations and $90$ seconds. On the capillary energy the subspace methods took $668$ to $863$ iterations and $67$ to $82$ seconds, L-BFGS $121$ seconds and the AGD variants $133$ to $201$ seconds. The subspace methods were thus the fastest in both tests.

With the rule of Proposition~\ref{thm:inexact} and $\gamma=1$, we initialize damped Newton at a gradient step whose length is adjusted by factors of two until the Armijo test with parameter $1/4$ holds and fails at twice the accepted length. Newton iterations continue until $\nu_k^2\le a_k\norm{\nabla f(x_k)}^2$. The activation schedule then needs $1$ to $1.5$ Newton steps per iteration on logistic regression ($\kappa=10^4,10^6$) and the capillary energy ($n=200$, $Q=40$), $5$ to $9$ times fewer than with exact subproblems, and never more outer iterations.

\section{Conclusion}\label{sec:conclusion}

We asked which subspaces built from the history contain the momentum of an accelerated method when $\kappa$ is unknown. For memories fixed in advance the answer is a time scale. They need a window of length $\Omega(\sqrt\kappa)$ or a forgetting index $\Omega(\kappa^{-1/2})$, which excludes SESOP, the plane form of conjugate gradients and GD\nobreakdash-ELS, also on smooth functions, and with one pair and a bounded window the scale is $1-\rho=\Theta(\kappa^{-1/2})$. This scale can be reached without $\kappa$: filter pairs created at fixed times are uniformly accelerated, and a single adaptive pair is accelerated up to a polynomial prefactor. We do not know whether SESOP is slow at moderate accuracy or whether an adapted constant memory can be uniformly accelerated. On the logistic-regression and capillary test problems the slow regime did not appear, as expected: the tails start after a $\kappa$-dependent reduction of the gap, and for logistic regression the local condition number is far below the nominal one.

\section*{Declarations}
\textbf{Funding.} The research was supported by Russian Science Foundation (project No.~\mbox{21-71-30005-$\Pi$}), \url{https://rscf.ru/en/project/21-71-30005/}. \textbf{Competing interests.} The authors declare no competing interests. \textbf{Data availability.} No external data were used. All test problems are generated by the code. \textbf{Code availability.} The code for all experiments, the computed results and instructions are available on Zenodo, \url{https://doi.org/10.5281/zenodo.23916709}. \textbf{Author contributions.} I.~Kuruzov and A.~Gasnikov conceived and formulated the problem. I.~Kuruzov obtained the main theoretical results, wrote the manuscript and prepared the presentation of the results. I.~Dorofeev carried out the numerical experiments and contributed to the auxiliary results. All authors read and approved the final manuscript. \textbf{Use of AI tools.} The authors verified every proof independently. An AI assistant was also used to write the code for the experiments, to recheck the proofs, and to help with writing. The authors are fully responsible for the paper.

{\hbadness=10000

}

\end{document}